\documentclass[11pt]{article}
\usepackage[margin=1in]{geometry}
\usepackage{amsfonts,amsmath,amssymb}
\usepackage{array}

\PassOptionsToPackage{obeyspaces}{url}
\usepackage[colorlinks=true,citecolor=blue,urlcolor=blue,linkcolor=blue,bookmarksopen=true]{hyperref}
\usepackage{amsthm}
\usepackage{tikz}
\usetikzlibrary{calc,shapes.geometric,positioning,decorations.pathmorphing}
\usepackage{mathtools}
\usepackage{subcaption}
\usepackage{enumitem}

\usepackage{realhats}

\usepackage{etoolbox}
\patchcmd{\thebibliography}{\leftmargin\labelwidth}{\leftmargin\labelwidth\addtolength\itemsep{-0.1\baselineskip}}{}{}

\author{%
    Joseph Briggs\thanks{Dept.\ of Mathematics \& Statistics, Auburn University, Auburn, AL, USA. \texttt{\{jgb0059,coc0014\}@auburn.edu}.}
    \and
    Alex Parker\thanks{Dept.\ of Mathematics, Iowa State University, Ames, IA, USA. \texttt{\{abparker,cschwi\}@iastate.edu}.} \thanks{Supported in part by U.S.\ taxpayers through NSF RTG Grant DMS-1839918.}
    \and
    Coy Schwieder\footnotemark[2] \footnotemark[3]
    \and
    Chris Wells\footnotemark[1] \footnotemark[3]
}

\title{Frogs, hats and common subsequences}
\date{}

\usepackage[nameinlink]{cleveref}

\newtheorem{theorem}{Theorem}

\newtheorem{lemma}[theorem]{Lemma}
\newtheorem{corollary}[theorem]{Corollary}

\newtheorem{proposition}[theorem]{Proposition}
\crefname{proposition}{proposition}{propositions}

\newtheorem{observation}[theorem]{Observation}
\crefname{observation}{observation}{observations}

\crefname{claim}{claim}{claims}

\crefname{conj}{conjecture}{conjectures}

\theoremstyle{definition}

\newtheorem{defn}[theorem]{Definition}
\crefname{defn}{definition}{definitions}

\crefname{remark}{remark}{remarks}

\newtheorem*{remark*}{Remark}

\Crefname{equation}{Eq.\!}{Eqs.\!}

\newlist{properties}{enumerate}{4}
\setlist[properties,1]{%
    label=(\arabic*),
    ref=\arabic*
}
\setlist[properties,2]{%
    label=(\roman*),
    ref=\thepropertiesi.\roman*
}
\crefname{propertiesi}{property}{properties}
\crefname{propertiesii}{property}{properties}
\crefname{propertiesiii}{property}{properties}
\crefname{propertiesiv}{property}{properties}

\newlist{rules}{enumerate}{2}
\setlist[rules,1]{label=\Alph*), ref=\Alph*}
\setlist[rules,2]{label=\therulesi\arabic*), ref=\therulesi\arabic*}
\crefname{rulesi}{rule}{rules}
\crefname{rulesii}{rule}{rules}

\newcommand{\crefOr}[1]{%
    \crefmultiformat{rulesii}{##2\csuse{cref@rulesii@name}~##1##3}%
    { or~##2##1##3}{, ##2##1##3}{ or~##2##1##3}%
    \cref{#1}%
    \crefmultiformat{rulesii}{\csuse{cref@rulesii@name@plural}~##2##1##3}%
    { and~##2##1##3}{, ##2##1##3}{, and~##2##1##3}%
}

\makeatother

\let\theparentequation\theequation
\patchcmd{\theparentequation}{equation}{parentequation}{}{}

\newcommand*{\eqdef}{\stackrel{\mbox{\normalfont\tiny def}}{=}}   
\newcommand*{\abs}[1]{\lvert #1\rvert}                           
\newcommand*{\Z}{\mathbb{Z}}                                     

\newcommand*{\oeis}[1]{\href{http://oeis.org/#1}{\texttt{OEIS:#1}}}

\DeclareMathOperator*{\E}{\mathbb{E}}                             
\DeclareMathOperator*{\Var}{Var}                                  
\DeclareMathOperator*{\len}{len}                                  
\DeclareMathOperator*{\LCS}{LCS}                                  
\def\mydot at (#1,#2){\fill (#1,#2) circle [radius=0.2]}         
\newcommand*{\bet}{\Sigma}										

\newcommand*{\refmapsto}[1]{\xmapsto{\text{\ref{#1}}}}

\let\oldhat\hat
\newcommand*{\crown}[1]{\oldhat[crown]{#1}\vphantom{#1}}
\renewcommand*{\hat}[1]{\oldhat[tophat]{#1}\vphantom{#1}}

\def\lilypad{\raisebox{-0.3ex}{\includegraphics[height=0.66em]{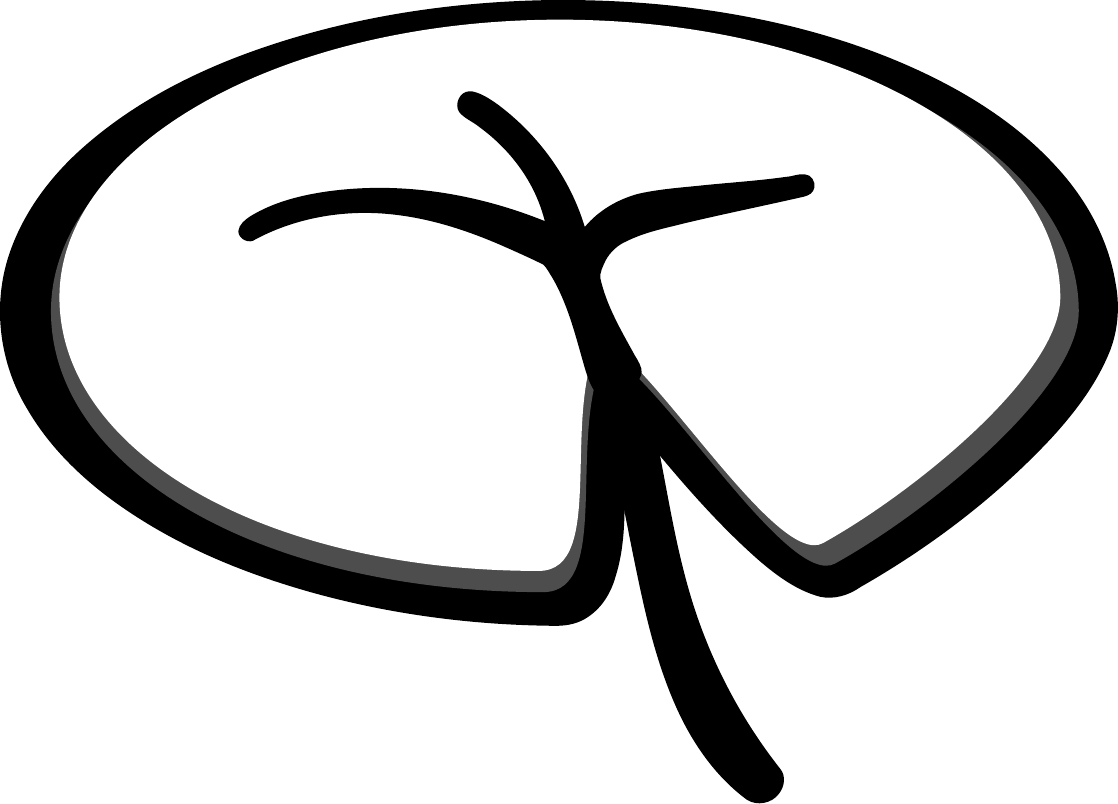}}}
\def\froggie{\protect\includegraphics[scale=0.014]{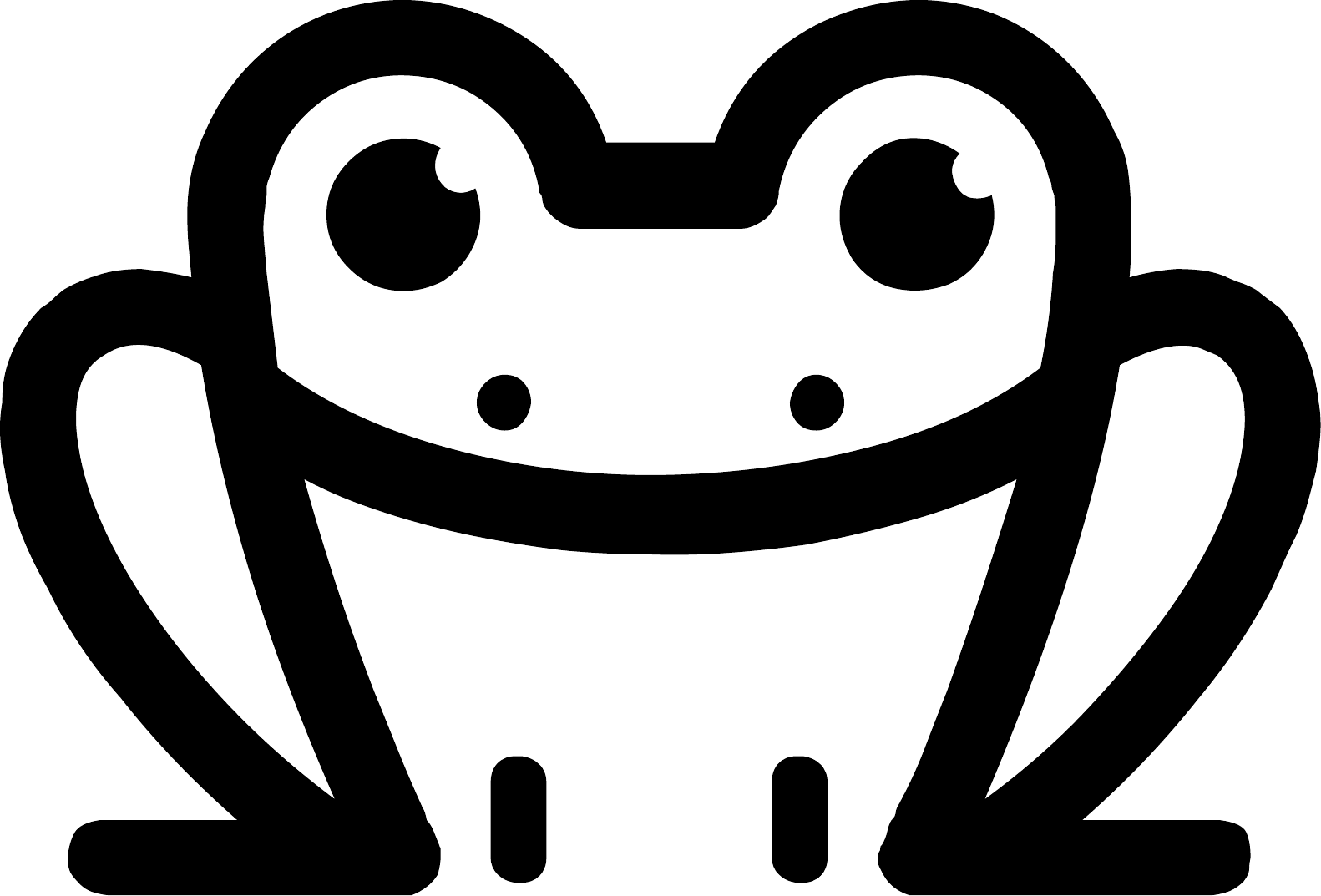}}
\def\froggiehat{\hat{\froggie}}
\def\froggiecrown{\crown{\froggie}}
\def\froggieb{\vphantom{\froggiehat}\froggie}

\def\nofroggie{\vphantom{\froggiehat}}

\newcommand*\Frogs{\mathcal{F}}
\newcommand*\blindFrogs[2]{\Frogs_{{#1},{#2}}}
\newcommand*\hattedFrogs[2]{\hat{\Frogs}_{{#1},{#2}}}
\newcommand*\crownedFrogs[2]{\crown{\Frogs}_{{#1},{#2}}}
\newcommand*\crownedFrogsSub[3]{\crownedFrogs{#2}{#3}^{\normalfont\texttt{#1}}}

\newcommand*\crownedFrogsStart[2]{\crownedFrogsSub{start}{#1}{#2}}
\newcommand*\crownedFrogsEnd[2]{\crownedFrogsSub{end}{#1}{#2}}

\newcommand*\poke[1]{\operatorname{\texttt{\textsc{Poke}}}_{#1}}
\newcommand*\dethrone{\operatorname{\texttt{\textsc{DeThrone}}}}
\newcommand*\move{\operatorname{\texttt{\textsc{Hop}}}}
\newcommand*\doff{\operatorname{\texttt{\textsc{Doff}}}}
\newcommand*\rot{\operatorname{\texttt{\textsc{Rotate}}}}
\newcommand*\agitated{{\normalfont\texttt{agitated}}}
\newcommand*\settled{{\normalfont\texttt{settled}}}
\DeclareMathOperator\disp{D}
\DeclareMathOperator\hop{H}
\DeclareMathOperator\Hop{\mathcal{H}}

\DeclareMathOperator\opp{opp}
\DeclareMathOperator\column{Col}
\newcommand*\eb[1]{\operatorname{e}^-_{#1}}

\newcommand*\sfcaption{\vskip-0.5cm}
\newcommand*\sfbetween{\vskip0.5cm}

\newcommand{\domino}[2]{%
    \scalebox{1}{
        \begin{tikzpicture}[baseline=(center)]
            \node (1) at (0,0.5) {$#1$};
            \node (2) at (0,0) {$#2$};
            \draw (-0.25,0.25)--(0.25,0.25);
            \draw (-0.25,-0.25)--(-0.25,0.75)--(0.25,0.75)--(0.25,-0.25)--(-0.25,-0.25);
            \coordinate (center) at (0,0.15);
        \end{tikzpicture}%
    }%
}
\newcommand{\dominoblock}[4]{%
    \scalebox{1}{
        \begin{tikzpicture}[baseline=(center)]
            \node (1) at (0,0.5) {$#1$};
            \node (2) at (0,0) {$#2$};
            \draw (-0.25,0.25)--(0.25,0.25);
            \draw (-0.25,-0.25)--(-0.25,0.75)--(0.25,0.75)--(0.25,-0.25)--(-0.25,-0.25);
            \node (3) at (0.5,0.5) {$#3$};
            \node (4) at (0.5,0) {$#4$};
            \draw (0.25,0.25)--(0.75,0.25);
            \draw (0.25,-0.25)--(0.25,0.75)--(0.75,0.75)--(0.75,-0.25)--(0.25,-0.25);
            \coordinate (center) at (0,0.15);
        \end{tikzpicture}%
    }%
}
\newcommand{\dominoTransition}[3]{%
    \scalebox{1}{
        \begin{tikzpicture}[baseline=(center)]
            \node (1) at (0.5,0.5) {$#1$};
            \node (2) at (0.5,0) {$#2$};
            \draw (0.25,0.25)--(0.75,0.25);
            \draw (0.25,-0.25)--(0.25,0.75)--(0.75,0.75)--(0.75,-0.25)--(0.25,-0.25);
            \node (3) at (0,1) {$#3$};
            \coordinate (center) at (0,0.15);
        \end{tikzpicture}%
    }%
}
\newcommand{\dominoTransitioned}[3]{%
    \scalebox{1}{
        \begin{tikzpicture}[baseline=(center)]
            \node (1) at (0.5,0.5) {$#1$};
            \node (2) at (0.5,0) {$#2$};
            \draw (0.25,0.25)--(0.75,0.25);
            \draw (0.25,-0.25)--(0.25,0.75)--(0.75,0.75)--(0.75,-0.25)--(0.25,-0.25);
            \node (4) at (0.5,1) {$#3$};
            \coordinate (center) at (0,0.15);
        \end{tikzpicture}%
    }%
}

\begin{document}

\maketitle

\begin{abstract}
    Write $W^{(n)}$ to mean the $n$-letter word obtained by repeating a fixed word $W$ and let $R_n$ denote a uniformly random $n$-letter word sampled from the same alphabet as $W$.
    We are interested in the average length of the longest common subsequence between $W^{(n)}$ and $R_n$, which is known to be $\gamma(W)\cdot n+o(n)$ for some constant $\gamma(W)$.
    Bukh and Cox recently developed an interacting particle system, dubbed the frog dynamics, which can be used to compute the constant $\gamma(W)$ for any fixed word $W$.
    They successfully analyzed the simplest case of the frog dynamics to find an explicit formula for the constants $\gamma(12\cdots k)$.
    We continue this study by using the frog dynamics to find an explicit formula for the constants $\gamma(12\cdots kk\cdots 21)$.
    The frog dynamics in this case is a variation of the PushTASEP on the ring where some clocks are identical.
    Interestingly, exclusion processes with correlated clocks of this type appear to have not been analyzed before.
    Our analysis leads to a seemingly new combinatorial object which could be of independent interest: frogs with hats!
\end{abstract}

\section{Introduction}

A \emph{word} is a finite sequence of letters from some alphabet.
A \emph{subsequence} of a word $W$ is a word obtained by deleting some letters from $W$; the letters in a subsequence are not required to appear contiguously in $W$.
A \emph{common subsequence} between words $W$ and $W'$ is a subsequence of both $W$ and $W'$.
We denote by $\LCS(W,W')$ the length of the \emph{longest common subsequence} between $W$ and $W'$.
The quantity $\LCS(W,W')$ is a measure of similarity between $W$ and $W'$ and is equivalent to the edit distance between $W$ and $W'$; that is, the fewest number of insertions/deletions needed to transform $W$ into $W'$.

Throughout the paper, we use $\bet$ to denote an alphabet (which is simply a finite set) and write $R\sim\bet^n$ to indicate that $R$ is a word chosen uniformly at random from $\bet^n$.

A fundamental, long-standing problem is to understand $\LCS(R,R')$ for a pair of independently sampled $R,R'\sim\bet^n$.
This problem was introduced by Chv\'atal and Sankoff~\cite{chvatal_1975} in 1975 and it is still poorly understood to date.
In \cite{chvatal_1975}, Chv\'atal and Sankoff proved that for any alphabet $\bet$, there exists a constant $\gamma=\gamma(\abs\bet)$ for which
\[
    \E_{R,R'\sim\bet^n}\LCS(R,R')=\gamma n+o(n).
\]
The constant $\gamma$ is now known as the Chv\'atal--Sankoff constant for $\bet$.
To date, not a single value of $\gamma$ is known.
Even in the case of the binary alphabet, the best known bounds are $0.788071\leq\gamma(2)\leq0.826280$, both due to Lueker~\cite{lueker_2009}.
However, Tiskin~\cite{tiskin_algebraic} recently announced that $\gamma(2)$ is a solution to an explicit (but large) system of polynomial equations.
While individual values of $\gamma$ remain unknown, Kiwi, Loebl and Matou\v{s}ek~\cite{kiwi_2005} established the asymptotic $\gamma(k)=\bigl(2+o(1)\bigr)/\sqrt{k}$ as $k\to\infty$.

Beyond the actual value of $\gamma$, the error-term is poorly understood as well.
The best bound is due to Alexander~\cite{alexander_1994} who established $\E\LCS(R,R')=\gamma n+O\bigl(\sqrt{n\log n}\bigr)$, but it is thought that perhaps $\E\LCS(R,R')=\gamma n-\Theta(n^{1/3})$ (see {\cite[Section 5]{bukh_frogs}}).

Additionally, it is thought that if $R,R'\sim\bet^n$, then $\LCS(R,R')$ should be approximately normal with linear variance; however, it is not even known if $\Var\LCS(R,R')$ tends toward infinity with $n$.
\medskip

The difficulty of understanding $\LCS(R,R')$ leads one to turn to understanding a class of related, but simpler, random variables: $\LCS(R,F)$ where $R$ is random and $F$ is ``fixed'' in some sense.
While one can fix the word $F$ in myriad ways, we will be concerned with the situation when $F$ is periodic.

For a word $W$, write $W^{(n)}$ to denote the periodic word of length $n$ which is obtained by repeating $W$ the appropriate number of times (which might be fractional if $\len W$ does not divide $n$).
For example, if $W=abba$, then $W^{(2)}=ab$ and $W^{(7)}=abbaabb$.

It is a routine exercise using superadditivity to show that for any alphabet $\bet$ and any finite word $W$, there is a constant $\gamma=\gamma(W,\bet)$ for which
\[
    \E_{R\sim\bet^n}\LCS(R,W^{(n)})=\gamma n+o(n).
\]
The leading constant $\gamma$ can be viewed as the Chv\'atal--Sankoff constant for the ``base-word'' $W$.
Unsurprisingly, actually computing the Chv\'atal--Sankoff constant and determining the order-of-magnitude of the error-term for any given base-word $W$ remains a highly non-trivial task.

The random variable $\LCS(R,W^{(n)})$ appears to have first been studied in its own right by Matzinger, Lember and Durringer~\cite{matzinger_variance}, who showed that it has linear variance when $R$ and $W$ come from the binary alphabet.
More recently, Bukh and Cox~\cite{bukh_frogs} developed an interacting particle system, dubbed the frog dynamics, which allows one to compute the Chv\'atal--Sankoff constant and the error term for any fixed base-word $W$.
Using the frog dynamics, they gave an explicit formula for the Chv\'atal--Sankoff constant and the order-of-magnitude of the error term in the case when $W=12\cdots k$ for some positive integer $k$.
Their results are not limited to the case when the random and periodic words have the same length; in the following, $\rho$ represents the ratio of their lengths.
\begin{theorem}[Bukh--Cox~\cite{bukh_frogs}]\label{bcLCS}
    Fix an integer $k\geq 2$ and a alphabet $\bet\supseteq[k]$.
    For any fixed $\rho\geq 0$,
    \[
        \E_{R\sim\bet^n}\LCS\bigl(R,(12\cdots k)^{(\rho n)}\bigr)=\biggl(\biggl(1-{m\over k}\biggr)\rho+{m\over \abs\bet (k+1-m)}\biggr)n-\tau\sqrt n+O(1),
    \]
    where $m\in[k]$ is the largest integer for which
    \begin{equation}\label{bc_threshold}
        {k(k+1)\over \abs\bet(k+2-m)(k+1-m)}\leq \rho,
    \end{equation}
    and $\tau=\tau(k,\abs\bet,\rho)\geq 0$ always with $\tau>0$ if and only if equality holds in \cref{bc_threshold}.
\end{theorem}

Prior to this current manuscript, the theorem above is the only instance of an infinite family of periodic words for which an explicit formula for their Chv\'atal--Sankoff constants is known.
This current paper is dedicated to expanding this list by understanding the case when $W=12\cdots kk\cdots 21$ for some positive integer $k$.

Before stating our main result, we remark that we attempted to understand many additional families of periodic words, but found success only with this particular one.
It is unclear if there is something particularly special about the base-word $W=12\cdots kk\cdots 21$ or if we simply lacked the proper insights to go further.

\begin{theorem}\label{lcsResult}
    Fix an integer $k\geq 2$ and a alphabet $\bet\supseteq[k]$.
    For any fixed $\rho\geq 0$,
    \[
        \E_{R\sim\bet^n}\LCS\bigl(R,(12\cdots kk\cdots 21)^{(\rho n)}\bigr)=\biggl(\biggl(1-{m\over 2k}\biggr)\rho +{\sum_{j\geq 0}{2k-2j\choose m-1-2j}\over \abs\bet\cdot\sum_{j\geq 0}{2k-2j\choose m-2j}}\biggr)n-\tau\sqrt n+O(1),
    \]
    where $m\in[2k]$ is the largest integer for which
    \begin{equation}\label{threshold}
        {2k\cdot\sum_{j\geq 0}{2k-2j\choose m-1-2j}\over \abs\bet\cdot\sum_{j\geq 0}{2k-2j\choose m-2j}}-{2k\cdot\sum_{j\geq 0}{2k-2j\choose m-2-2j}\over \abs\bet\cdot\sum_{j\geq 0}{2k-2j\choose m-1-2j}}\leq \rho,
    \end{equation}
    and $\tau=\tau(k,\abs\bet,\rho)\geq 0$ always with $\tau>0$ if and only if equality holds in \cref{threshold}.
\end{theorem}

This paper is organized as follows.
In \Cref{sec:preliminaries}, we outline the Bukh--Cox frog dynamics and demonstrate how it will be used to prove \Cref{lcsResult}.
Then, in \Cref{sec:hattedFrogs}, we introduce a seemingly new combinatorial object, which may be of independent interest, in order to analyze the frog dynamics: frogs with hats!
Finally, in \Cref{sec:speeds}, we use the frogs with hats to compute the necessary quantities of the frog dynamics which yield \Cref{lcsResult}.

\subsection{Notation}
For a probability measure $\mu$, we write $x\sim\mu$ to indicate that $x$ is distributed according to $\mu$.
For a finite set $\Omega$, we write $x\sim\Omega$ to indicate that $x$ is drawn uniformly at random from $\Omega$.

For a set $\Omega$ and an element $x$, we write $\Omega+x$ to indicate $\Omega\cup\{x\}$ and write $\Omega-x$ to indicate $\Omega\setminus\{x\}$.
We often combine these two operations, e.g.\ $\Omega-x+y=(\Omega\setminus\{x\})\cup\{y\}$.
These abbreviations are used to aid readability and will be avoided if they ever introduce ambiguity.

When two sets $\Omega_1,\Omega_2$ are disjoint, we will often write their union as $\Omega_1\sqcup\Omega_2$ to emphasize the fact that the two sets are disjoint.

For integers $m,n$, we use the shorthand $[m,n]\eqdef\{m,m+1,\dots,n-1,n\}$ where $[m,n]=\varnothing$ if $n<m$.
As is common, we abbreviate $[n]\eqdef[1,n]=\{1,2,\dots,n\}$.
For a set $\Omega$ and an integer $k$, we denote by ${\Omega\choose k}$ the set of all $k$-sized subsets of $\Omega$.
We additionally denote the set of all subsets of size at most $k$ by ${\Omega\choose\leq k}$.

\section{Preliminaries}\label{sec:preliminaries}
We begin by outlining the Bukh--Cox frog dynamics from \cite{bukh_frogs}.

Fix a word $W\in\bet^\ell$ where $\bet$ is some alphabet.
Imagine a circle of $\ell$ lily pads, $\lilypad_0,\dots,\lilypad_{\ell-1}$, arranged in clockwise order, where $\lilypad_i$ is labeled with the $(i+1)$th letter of $W$.
Each of these lily pads is additionally occupied by a frog.
The $\ell$ frogs vary from a large nasty frog to a little harmless froggie.
No two frogs are equally nasty and so we denote them by $\froggie_1,\dots,\froggie_\ell$ where $\froggie_1$ is the nastiest frog and $\froggie_\ell$ is the least nasty frog.
The arrangement of the frogs on the lily pads is described by a bijection $F\colon\{\froggie_1,\dots,\froggie_\ell\}\to\{\lilypad_0,\dots,\lilypad_{\ell-1}\}$; any such bijection is called a \emph{frog arrangement}.

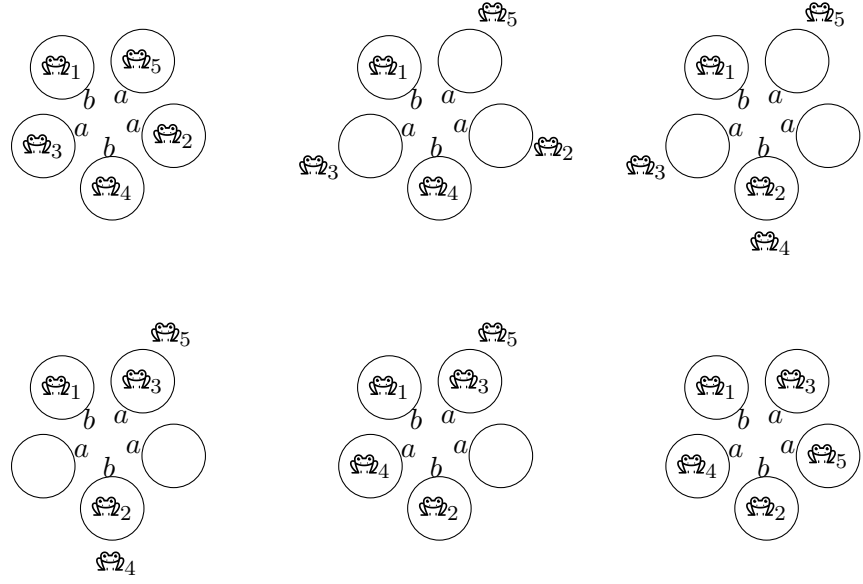
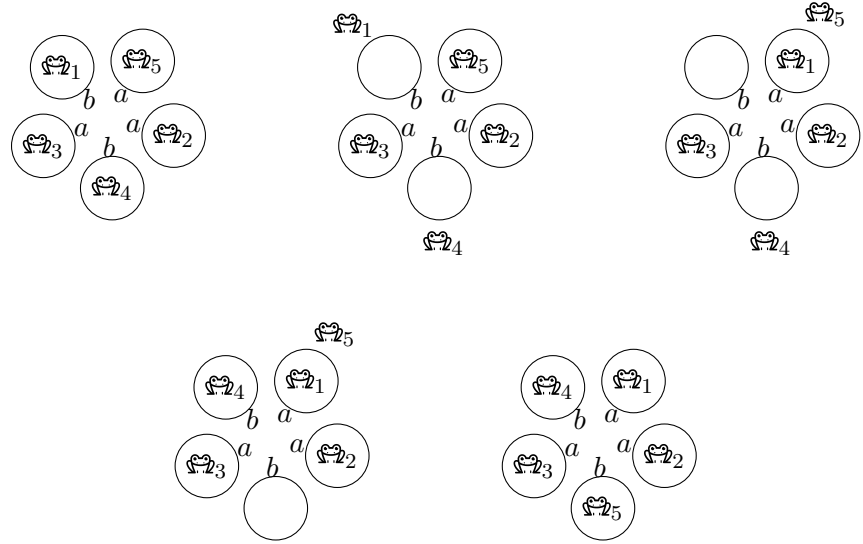
\begin{figure}[tp]
    \def\labelrad{0.4}
    \begin{subfigure}{\textwidth}
        \begin{center}
            \begin{tikzpicture}[scale = 0.7]
                \clip (-3,-3) rectangle (3,3);
                \foreach \x in {1,2,3,4,5} {
                    \coordinate (0\x) at ({180+45/2-\x*360/5}:1.3);
                    \draw (0\x) circle (0.6);
                }
                \node (l1) at ($\labelrad*(01)$) {$b$};
                \node (l2) at ($\labelrad*(02)$) {$a$};
                \node (l3) at ($\labelrad*(03)$) {$a$};
                \node (l4) at ($\labelrad*(04)$) {$b$};
                \node (l5) at ($\labelrad*(05)$) {$a$};
                \node (f1) at (01) {$\froggie_1$};
                \node (f2) at (02) {$\froggie_5$};
                \node (f3) at (03) {$\froggie_2$};
                \node (f4) at (04) {$\froggie_4$};
                \node (f5) at (05) {$\froggie_3$};
            \end{tikzpicture}
        \begin{tikzpicture}[scale = 0.7]
            \clip (-3,-3) rectangle (3,3);
            \foreach \x in {1,2,3,4,5} {
                \coordinate (0\x) at ({180+45/2-\x*360/5}:1.3);
                \draw (0\x) circle (0.6);
            }
            \node (l1) at ($\labelrad*(01)$) {$b$};
            \node (l2) at ($\labelrad*(02)$) {$a$};
            \node (l3) at ($\labelrad*(03)$) {$a$};
            \node (l4) at ($\labelrad*(04)$) {$b$};
            \node (l5) at ($\labelrad*(05)$) {$a$};
            \node (f1) at (01) {$\froggie_1$};
            \node (f2) at ($1.8*(02)$) {$\froggie_5$};
            \node (f3) at ($1.8*(03)$) {$\froggie_2$};
            \node (f4) at (04) {$\froggie_4$};
            \node (f5) at ($1.8*(05)$) {$\froggie_3$};
        \end{tikzpicture}
        \begin{tikzpicture}[scale = 0.7]
            \clip (-3,-3) rectangle (3,3);
            \foreach \x in {1,2,3,4,5} {
                \coordinate (0\x) at ({180+45/2-\x*360/5}:1.3);
                \draw (0\x) circle (0.6);
            }
            \node (l1) at ($\labelrad*(01)$) {$b$};
            \node (l2) at ($\labelrad*(02)$) {$a$};
            \node (l3) at ($\labelrad*(03)$) {$a$};
            \node (l4) at ($\labelrad*(04)$) {$b$};
            \node (l5) at ($\labelrad*(05)$) {$a$};
            \node (f1) at (01) {$\froggie_1$};
            \node (f2) at ($1.8*(02)$) {$\froggie_5$};
            \node (f3) at (04) {$\froggie_2$};
            \node (f4) at ($1.8*(04)$) {$\froggie_4$};
            \node (f5) at ($1.8*(05)$) {$\froggie_3$};
        \end{tikzpicture}
        \begin{tikzpicture}[scale = 0.7]
            \clip (-3,-3) rectangle (3,3);
            \foreach \x in {1,2,3,4,5} {
                \coordinate (0\x) at ({180+45/2-\x*360/5}:1.3);
                \draw (0\x) circle (0.6);
            }
            \node (l1) at ($\labelrad*(01)$) {$b$};
            \node (l2) at ($\labelrad*(02)$) {$a$};
            \node (l3) at ($\labelrad*(03)$) {$a$};
            \node (l4) at ($\labelrad*(04)$) {$b$};
            \node (l5) at ($\labelrad*(05)$) {$a$};
            \node (f1) at (01) {$\froggie_1$};
            \node (f2) at ($1.8*(02)$) {$\froggie_5$};
            \node (f3) at (04) {$\froggie_2$};
            \node (f4) at ($1.8*(04)$) {$\froggie_4$};
            \node (f5) at (02) {$\froggie_3$};
        \end{tikzpicture}
        \begin{tikzpicture}[scale = 0.7]
            \clip (-3,-3) rectangle (3,3);
            \foreach \x in {1,2,3,4,5} {
                \coordinate (0\x) at ({180+45/2-\x*360/5}:1.3);
                \draw (0\x) circle (0.6);
            }
            \node (l1) at ($\labelrad*(01)$) {$b$};
            \node (l2) at ($\labelrad*(02)$) {$a$};
            \node (l3) at ($\labelrad*(03)$) {$a$};
            \node (l4) at ($\labelrad*(04)$) {$b$};
            \node (l5) at ($\labelrad*(05)$) {$a$};
            \node (f1) at (01) {$\froggie_1$};
            \node (f2) at ($1.8*(02)$) {$\froggie_5$};
            \node (f3) at (04) {$\froggie_2$};
            \node (f4) at (05) {$\froggie_4$};
            \node (f5) at (02) {$\froggie_3$};
        \end{tikzpicture}
        \begin{tikzpicture}[scale = 0.7]
            \clip (-3,-3) rectangle (3,3);
            \foreach \x in {1,2,3,4,5} {
                \coordinate (0\x) at ({180+45/2-\x*360/5}:1.3);
                \draw (0\x) circle (0.6);
            }
            \node (l1) at ($\labelrad*(01)$) {$b$};
            \node (l2) at ($\labelrad*(02)$) {$a$};
            \node (l3) at ($\labelrad*(03)$) {$a$};
            \node (l4) at ($\labelrad*(04)$) {$b$};
            \node (l5) at ($\labelrad*(05)$) {$a$};
            \node (f1) at (01) {$\froggie_1$};
            \node (f2) at (03) {$\froggie_5$};
            \node (f3) at (04) {$\froggie_2$};
            \node (f4) at (05) {$\froggie_4$};
            \node (f5) at (02) {$\froggie_3$};
        \end{tikzpicture}
        \end{center}
        \sfcaption
        \caption{
            The monster pokes lily pads labeled $a$!
            The top-left image is $F$ and the bottom-right image is $Fa$.
            Here, $\disp_1(F,a)=0$, $\disp_2(F,a)=1$, $\disp_3(F,a)=2$, $\disp_4(F,a)=1$ and $\disp_5(F,a)=1$.
        }
    \end{subfigure}
    \sfbetween

    \begin{subfigure}{\textwidth}
        \begin{center}
            \begin{tikzpicture}[scale = 0.7]
                \clip (-3,-3) rectangle (3,3);
                \foreach \x in {1,2,3,4,5} {
                    \coordinate (0\x) at ({180+45/2-\x*360/5}:1.3);
                    \draw (0\x) circle (0.6);
                }
                \node (l1) at ($\labelrad*(01)$) {$b$};
                \node (l2) at ($\labelrad*(02)$) {$a$};
                \node (l3) at ($\labelrad*(03)$) {$a$};
                \node (l4) at ($\labelrad*(04)$) {$b$};
                \node (l5) at ($\labelrad*(05)$) {$a$};
                \node (f1) at (01) {$\froggie_1$};
                \node (f2) at (02) {$\froggie_5$};
                \node (f3) at (03) {$\froggie_2$};
                \node (f4) at (04) {$\froggie_4$};
                \node (f5) at (05) {$\froggie_3$};
            \end{tikzpicture}
            \begin{tikzpicture}[scale = 0.7]
                \clip (-3,-3) rectangle (3,3);
                \foreach \x in {1,2,3,4,5} {
                    \coordinate (0\x) at ({180+45/2-\x*360/5}:1.3);
                    \draw (0\x) circle (0.6);
                }
                \node (l1) at ($\labelrad*(01)$) {$b$};
                \node (l2) at ($\labelrad*(02)$) {$a$};
                \node (l3) at ($\labelrad*(03)$) {$a$};
                \node (l4) at ($\labelrad*(04)$) {$b$};
                \node (l5) at ($\labelrad*(05)$) {$a$};
                \node (f1) at ($1.8*(01)$) {$\froggie_1$};
                \node (f2) at (02) {$\froggie_5$};
                \node (f3) at (03) {$\froggie_2$};
                \node (f4) at ($1.8*(04)$) {$\froggie_4$};
                \node (f5) at (05) {$\froggie_3$};
            \end{tikzpicture}
            \begin{tikzpicture}[scale = 0.7]
                \clip (-3,-3) rectangle (3,3);
                \foreach \x in {1,2,3,4,5} {
                    \coordinate (0\x) at ({180+45/2-\x*360/5}:1.3);
                    \draw (0\x) circle (0.6);
                }
                \node (l1) at ($\labelrad*(01)$) {$b$};
                \node (l2) at ($\labelrad*(02)$) {$a$};
                \node (l3) at ($\labelrad*(03)$) {$a$};
                \node (l4) at ($\labelrad*(04)$) {$b$};
                \node (l5) at ($\labelrad*(05)$) {$a$};
                \node (f1) at (02) {$\froggie_1$};
                \node (f2) at ($1.8*(02)$) {$\froggie_5$};
                \node (f3) at (03) {$\froggie_2$};
                \node (f4) at ($1.8*(04)$) {$\froggie_4$};
                \node (f5) at (05) {$\froggie_3$};
            \end{tikzpicture}
            \begin{tikzpicture}[scale = 0.7]
                \clip (-3,-3) rectangle (3,3);
                \foreach \x in {1,2,3,4,5} {
                    \coordinate (0\x) at ({180+45/2-\x*360/5}:1.3);
                    \draw (0\x) circle (0.6);
                }
                \node (l1) at ($\labelrad*(01)$) {$b$};
                \node (l2) at ($\labelrad*(02)$) {$a$};
                \node (l3) at ($\labelrad*(03)$) {$a$};
                \node (l4) at ($\labelrad*(04)$) {$b$};
                \node (l5) at ($\labelrad*(05)$) {$a$};
                \node (f1) at (02) {$\froggie_1$};
                \node (f2) at ($1.8*(02)$) {$\froggie_5$};
                \node (f3) at (03) {$\froggie_2$};
                \node (f4) at (01) {$\froggie_4$};
                \node (f5) at (05) {$\froggie_3$};
            \end{tikzpicture}
            \begin{tikzpicture}[scale = 0.7]
                \clip (-3,-3) rectangle (3,3);
                \foreach \x in {1,2,3,4,5} {
                    \coordinate (0\x) at ({180+45/2-\x*360/5}:1.3);
                    \draw (0\x) circle (0.6);
                }
                \node (l1) at ($\labelrad*(01)$) {$b$};
                \node (l2) at ($\labelrad*(02)$) {$a$};
                \node (l3) at ($\labelrad*(03)$) {$a$};
                \node (l4) at ($\labelrad*(04)$) {$b$};
                \node (l5) at ($\labelrad*(05)$) {$a$};
                \node (f1) at (02) {$\froggie_1$};
                \node (f2) at (04) {$\froggie_5$};
                \node (f3) at (03) {$\froggie_2$};
                \node (f4) at (01) {$\froggie_4$};
                \node (f5) at (05) {$\froggie_3$};
            \end{tikzpicture}
        \end{center}
        \sfcaption
        \caption{
            The monster pokes lily pads labeled $b$!
            The top-left image is $F$ and the bottom-right image is $Fb$.
            Here, $\disp_1(F,b)=1$, $\disp_2(F,b)=0$, $\disp_3(F,b)=0$, $\disp_4(F,b)=2$ and $\disp_5(F,b)=2$.
        }
    \end{subfigure}
    \caption{\label{fig:frogProcess}
        Examples of the frog process.
        An agitated frog is drawn just outside of the lily pad it occupies.
        The frogs hop in the clockwise direction.
    }
\end{figure}
Given a frog arrangement, the \emph{frog process} evolves as follows:
\begin{enumerate}
    \item The monster which lives in the pond selects some letter $a\in\bet$ and pokes all lily pads labeled $a$.
        This agitates any frog residing on a lily pad labeled $a$ and makes it want to jump away.
    \item In order of descending nastiness, starting from the nastiest frog $\froggie_1$, each of the agitated frogs will leap to the next `available' lily pad; that is, the closest lily pad in the clockwise order that is either empty or occupied by a less nasty frog.
        Upon leaping to this next available lily pad, the frog that just hopped calms down and any previous occupant becomes agitated.

        This repeats until all frogs are calm once more.
\end{enumerate}
Note that, with each step of the process, the nastiest agitated frog gets less and less nasty.
This guarantees the termination of the process and that no frog jumps over another agitated frog.
\Cref{fig:frogProcess} includes two examples of the frog process.

For a frog arrangement $F$ and a letter $a$, let $Fa$ denote the frog arrangement that results from starting with $F$, poking all lily pads labeled $a$ and waiting for the ensuing frenzy to settle.
Note that if no lily pads are labeled $a$, then $Fa=F$.
For each $m\in[\ell]$, let $\disp_m(F,a)$ denote the total number of lily pads that $\froggie_m$ hopped in the transition from $F$ to $Fa$.
For example, $\disp_1(F,a)$ is $1$ if $F(\froggie_1)$ is labeled $a$ and is otherwise $0$.

For a word $R=r_1\cdots r_n\in\bet^n$, we write $FR\eqdef (((Fr_1)r_2)\cdots) r_n$.
That is $FR$ is the frog arrangement which results from first poking the letter $r_1$, then poking $r_2$, etc.
Note that this notation respects concatenation: $(FR)R'=F(RR')$ for any two words $R,R'$.
We additionally set $\disp_m(F,R)$ to be the total displacement of $\froggie_m$ upon poking $r_1$, then poking $r_2$, etc.
Note that $\disp_m(F,RR')=\disp_m(F,R)+\disp_m(FR,R')$ for any two words $R,R'$.

Starting with a frog arrangement $F_0$, set $F_n=F_0r_1r_2\cdots r_n$ where the letters $r_1,r_2,\dots$ are chosen independently and uniformly from $\bet$.
Since $F_n=F_{n-1}r_n$, the sequence $F_0,F_1,F_2,\ldots$ forms a Markov chain, which is known as the \emph{$(W,\bet)$ frog dynamics}.
\medskip

\begin{remark*}
The frog dynamics is a particular variation on the so-called ``colored/multi-species PushTASEP on the ring''.
We direct the reader to \cite{aggarwal_asep,ayyer_tasep} for more information about this general framework, but we do not assume any prior knowledge of this subject in this paper.
For those familiar with exclusion processes, though, we remark that this current paper appears to be the first time that such a process having ``correlated clocks'' have been considered.
\end{remark*}

A word $W$ is said to be \emph{irreducible} if it cannot be expressed as $W=UU\cdots U$ for some other word $U$.
Formally, $W$ is irreducible if there is \emph{no} word $U\neq W$ with $\len U\mid\len W$ for which $W=U^{(\len W)}$.
When discussing periodic words, we may always suppose that the base-word is irreducible since if $W=UU\cdots U$, then $W^{(n)}=U^{(n)}$ for every $n$.

Bukh and Cox established the following connection between longest common subsequences and the frog dynamics.
\begin{theorem}[Bukh--Cox~\cite{bukh_frogs}]\label{frogTheorem}
    Let $\bet$ be a alphabet and fix an irreducible word $W\in\bet^\ell$.
    \begin{enumerate}
        \item The $(W,\bet)$ frog dynamics admits a unique stationary distribution $\pi$.
        \item The average speed of $\froggie_m$, defined by
            \[
                s_m\eqdef\lim_{n\to\infty}\E_{R\sim\bet^n}{\disp_m(F_0,R)\over n},
            \]
            exists and is independent of the initial frog arrangement $F_0$.
            Furthermore,
            \[
                s_m=\E_{a\sim\bet}\ \E_{F\sim\pi}\disp_m(F,a)\qquad\text{for each }m\in[\ell].
            \]
        \item For every fixed $\rho\geq 0$,
            \[
                \E_{R\sim\bet^n}\LCS(R,W^{(\rho n)})=\biggl(\rho -{1\over\ell}\sum_{s_m\leq\rho}(\rho-s_m)\biggr)n-\tau\sqrt{n}+O(1)
            \]
            where $\tau=\tau(W,\bet,\rho)\geq 0$ always and $\tau>0$ if and only if $\rho\in\{s_1,\dots,s_\ell\}$.
    \end{enumerate}
\end{theorem}

Therefore, in order to determine the Chv\'atal--Sankoff constant for the base-word $W$, one needs only to compute the speeds $s_1,\dots,s_{\len W}$.

Bukh and Cox established \Cref{bcLCS} by showing that $s_m={k(k+1)\over\abs\bet(k+2-m)(k+1-m)}$ when $W=12\cdots k$ and $\bet\supseteq[k]$.
We will likewise establish \Cref{lcsResult} by showing that
\[
    s_m={2k\cdot\sum_{j\geq 0}{2k-2j\choose m-1-2j}\over \abs\bet\cdot\sum_{j\geq 0}{2k-2j\choose m-2j}}-{2k\cdot\sum_{j\geq 0}{2k-2j\choose m-2-2j}\over \abs\bet\cdot\sum_{j\geq 0}{2k-2j\choose m-1-2j}},
\]
when $W=12\cdots kk\cdots 21$ and $\bet\supseteq[k]$.
This is the content of \Cref{cumulative}.

\subsection{The blind-frog process and dynamics}

Fix a word $W\in\bet^\ell$ and an integer $m\in[\ell]$ and consider the frog dynamics associated with $W$.
One key idea which greatly simplifies the computations involved in analyzing the frog dynamics is a coupling with a collection of simpler chains, which we call the \emph{$m$-blind-frog dynamics}.
In the $m$-blind-frog dynamics, we focus only on $\froggie_{\leq m}\eqdef\{\froggie_1,\dots,\froggie_m\}$ , suppressing the distinction between $\froggie_1$ through $\froggie_m$.
As we are not interested in the position of any individual frog, we record the state of the system as a set $S\in{\{\lilypad_0,\dots,\lilypad_{\ell-1}\}\choose m}$.
We call such a set $S$ an \emph{$m$-blind-frog arrangement}, and we continue to refer to the elements of $S$ as `frogs', despite not knowing their relative nastiness.

Given an $m$-blind-frog arrangement $S\in{\{\lilypad_0,\dots,\lilypad_{\ell-1}\}\choose m}$, the \emph{$m$-blind-frog process} evolves as follows:
\begin{enumerate}
    \item The monster which lives in the pond selects some letter $a\in\bet$ and pokes all lily pads labeled $a$.
        This agitates any frog residing on a lily pad labeled $a$ and makes it want to jump away.
    \item With no restrictions on priority, an agitated frog will hop to the lily pad immediately succeeding it in the clockwise order and calm down.
        If this lily pad was previously occupied, then the previous occupant becomes agitated.

        This repeats until all $m$ frogs are content once more.
\end{enumerate}

For an $m$-blind-frog arrangement $S$ and $a\in\bet$, let $Sa$ be the $m$-blind-frog arrangement that results from starting with $S$, poking all lily pads labeled $a$ and waiting for the ensuing frenzy to settle.
Additionally, define $\hop(S,a)$ to be the total number of frogs that hopped in the transition from $S$ to $Sa$.

\begin{proposition}\label{blind_order}
    For any word $W\in\bet^\ell$ and any $m\in[\ell]$, the $m$-blind-frog process is well-defined.
    That is, for each $S\in{\{\lilypad_0,\dots,\lilypad_{\ell-1}\}\choose m}$ and each $a\in\bet$, the order in which the agitated frogs hop once the lily pads labeled $a$ are poked does not influence $Sa$ nor $\hop(S,a)$.
\end{proposition}
\begin{proof}
    Let $H\subseteq S$ be the union of all cyclic intervals contained in $S$ starting at a lily pad labeled $a$.
    Formally, $\lilypad_x\in H$ if and only if there is some $y\in\Z/\ell\Z$ (possibly $y=x$) for which $W_y=a$ and $\lilypad_y,\lilypad_{y+1},\dots,\lilypad_x\in S$.
    The effect of poking all lily pads labeled $a$ and waiting for the blind frogs to settle is precisely that every frog on a lily pad in $H$ hops exactly one lily pad forward, regardless of the order in which they hop.
    Indeed, every frog on a lily pad in $H$ will eventually become agitated, and no other frogs can become agitated since, by definition, there exists an empty lily pad between any element of $H$ and any element of $S\setminus H$ in the cyclic order.
    This observation implies that $Sa=(S\setminus H)\cup\{\lilypad_{x+1}:\lilypad_x\in H\}$ and $\hop(S,a)=\abs H$.
\end{proof}

Since the order in which the blind frogs hop is irrelevant, the blind-frog process aligns with the frog process.
We remark that this type of observation is nothing new and some version of it is likely well-known to anyone who studies related interacting particle systems: {\cite[Lemma 16]{ayyer_tasep}} gives an overview of the standard arguments.
See also {\cite[Appendix A]{aggarwal_asep}} for a more-involved version of this ``color merging'' idea.

\begin{proposition}\label{blind_frogs}
    If $W\in\bet^\ell$, then for any frog arrangement $F\colon \{\froggie_1,\dots,\froggie_\ell\}\to\{\lilypad_0,\dots,\lilypad_{\ell-1}\}$, $m\in[\ell]$ and $a\in\bet$, we have
    \begin{align*}
        (Fa)(\froggie_{\leq m}) &=\bigl(F(\froggie_{\leq m})\bigr)a,\qquad\text{and}\\
        \sum_{i=1}^m \disp_i(F,a) &= \hop(F(\froggie_{\leq m}),a).
    \end{align*}
\end{proposition}

\begin{proof}
    We first observe that the movement of $\froggie_1,\dots,\froggie_m$ is unaffected by the floundering of the less nasty frogs $\froggie_{m+1},\dots,\froggie_k$.
    Indeed, for any $i<j$, $\froggie_j$ can never agitate $\froggie_i$ and when $\froggie_i$ looks to hop, it never considers the presence of $\froggie_j$.
    Therefore, since we seek to track only the nastiest $m$ frogs, we may pretend as if $\froggie_{m+1},\dots,\froggie_k$ do not exist.

    Now, consider a slight modification to the frog process where, when $\froggie_i$ is agitated, instead of $\froggie_i$ hopping over any nastier frog, it instead hops onto the very next lily pad.
    If that lily pad is empty, $\froggie_i$ stops; otherwise there is another frog $\froggie_j$ occupying that lily pad.
    If $i>j$, then $\froggie_i$ settles down and $\froggie_j$ becomes agitated.
    Otherwise, $\froggie_i$ remains agitated and will continue to hop.

    This alternative viewpoint is readily seen to be equivalent to the original frog dynamics: we are simply considering a long leap to be comprised of smaller hops.
    Furthermore, with this alternative viewpoint, whenever a frog lands on an occupied lily pad, one of these two frogs hops away on the next step.
    Suppressing the distinction between $\froggie_1,\dots,\froggie_m$ does not change this fact and so, since the order in which the blind frogs hop is irrelevant in the blind-frog process (\Cref{blind_order}), the conclusion follows.
\end{proof}

Starting with an $m$-blind-frog arrangement $F_0$, set $F_n=F_0r_1r_2\cdots r_n$ where the letters $r_1,r_2,\ldots$ are chosen independently and uniformly from $\bet$.
Then the sequence $F_0,F_1,F_2,\ldots$ forms a Markov chain, which we call the \emph{$(W,\bet)$ $m$-blind-frog dynamics}.

Since the frog process aligns with the $m$-blind-frog process, we can couple the $(W,\bet)$ frog dynamics with the $(W,\bet)$ $m$-blind-frog dynamics by simply poking the same letter in each chain at each time step.

\begin{corollary}\label{speeds_v_hops}
    Fix an alphabet $\bet$ and an irreducible word $W\in\bet^\ell$
    If $\pi$ denotes the stationary distribution of the $(W,\bet)$ frog dynamics and $\pi_m$ denotes the stationary distribution of the $(W,\bet)$ $m$-blind-frog dynamics, then
    \begin{enumerate}
        \item $F\sim\pi\implies F(\froggie_{\leq m})\sim\pi_m$, and
        \item $\displaystyle \sum_{i=1}^m s_i=\E_{a\sim\bet}\ \E_{F\sim\pi_m}\hop(F,a)$.
    \end{enumerate}
\end{corollary}
In the above statement, note that since $\pi$ exists and is unique (\Cref{frogTheorem}), the same is true of $\pi_m$.

In the special case where $W=12\cdots k$ and $\Sigma\supseteq[k]$, Bukh and Cox found that the stationary distribution of the $m$-blind-frog dynamics is uniform~{\cite[Theorem 30]{bukh_frogs}}.
This fact allowed them to compute $\sum_{i=1}^m s_i$ for each $m\in[k]$ without much additional effort.
Alas, the stationary distribution of the $m$-blind-frog dynamics for $W=12\cdots kk\cdots 21$ is not anywhere as easy to describe, let alone prove.

\subsection{The blind-frogs associated with \texorpdfstring{$W=12\cdots kk\cdots 21$}{W=12...kk...21}}\label{ring_to_grid}

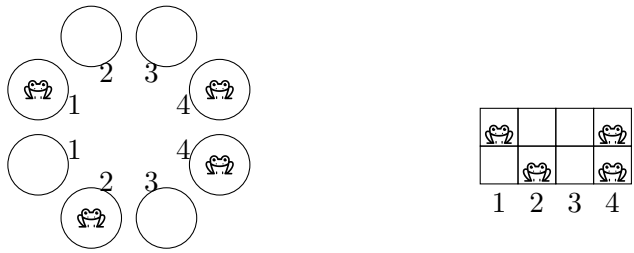
\begin{figure}[ht]
    \begin{center}
        \begin{tikzpicture}
            \clip (-2,-2) rectangle (2,2);
            \foreach \x in {1,2,3,4} {
                \coordinate (0\x) at ({180+45/2-\x*360/8}:1.3);
                \draw (0\x) circle (0.4);
                \node (l\x) at ($0.6*(0\x)$) {$\x$};
                \coordinate (1\x) at ({180-45/2+\x*360/8}:1.3);
                \draw (1\x) circle (0.4);
                \node (r\x) at ($0.6*(1\x)$) {$\x$};
            }
            \node (f1) at (01) {$\froggie$};
            \node (f2) at (04) {$\froggie$};
            \node (f3) at (14) {$\froggie$};
            \node (f4) at (12) {$\froggie$};
        \end{tikzpicture}
        \qquad
        \begin{tikzpicture}
            \clip (-2,-1.5) rectangle (2,2.5);
            \foreach \x in {0,0.5,1,1.5} {
                \draw ({\x-0.25},-0.25) rectangle (\x+0.25,0.25);
                \draw ({\x-0.25},0.25) rectangle (\x+0.25,0.75);
            }
            \foreach \x in {1,2,3,4} {
                \node (\x) at ({(\x-1)/2},-0.5) {$\x$};
            }
            \node (f1) at (0,0.5) {$\froggieb$};
            \node (f2) at (1.5,0.5) {$\froggieb$};
            \node (f3) at (1.5,0) {$\froggieb$};
            \node (f4) at (0.5,0) {$\froggieb$};
        \end{tikzpicture}
    \end{center}
    \caption{\label{squish_ring} Representing the ring of lily pads (left) by a grid (right).}
\end{figure}

When working with the word $W=12\cdots kk\cdots 21$, it is convenient to represent the ring of lily pads instead as the grid $[2]\times[k]$ (see \Cref{squish_ring}).

When considering frogs arranged on $[2]\times[k]$, the frogs in the top row (squares indexed by $(1,i)$) move to the right and frogs in the bottom row (squares indexed by $(2,i)$) move to the left, wrapping around at the sides.
In this way, the letter $c\in[k]$ corresponds directly to column $c$ in the grid.
Thus, instead of poking individual lily pads, we think about poking an entire column in this grid.

We denote the set of all arrangements of $m$ (blind) frogs on the grid $[2]\times[k]$ by $\blindFrogs km$; formally $\blindFrogs km\eqdef{[2]\times[k]\choose m}$.
\medskip

The most important observation about this grid arrangement is its symmetry.
Consider the map $\rot\colon[2]\times[k]\to[2]\times[k]$ defined by $\rot(i,j)=\rot(3-i, k+1-j)$, thereby rotating the grid by $180^\circ$.
The map $\rot$ can naturally be applied to each $F\in\blindFrogs km$ by simply applying $\rot$ pointwise to $F$.

\begin{observation}\label{rot_symmetry}
    Fix a letter $c\in[k]$ and set $c'=k+1-c$.
    Then for any $F\in\blindFrogs km$, we have $\rot(F)\in\blindFrogs km$ and $\rot(Fc)=\rot(F)c'$.
\end{observation}

\section{Frogs...with hats!}\label{sec:hattedFrogs}

From here on, whenever we refer to a `frog', that frog is blind in the sense that we do not care about relative nastiness.
Recall from \Cref{ring_to_grid} that we consider these frogs to live on the grid $[2]\times[k]$ instead of on a ring of lily pads.
In order to analyze the $m$-blind-frog dynamics associated with the word $W=12\cdots kk\cdots 21$, we need to introduce a new combinatorial object: frogs with hats!

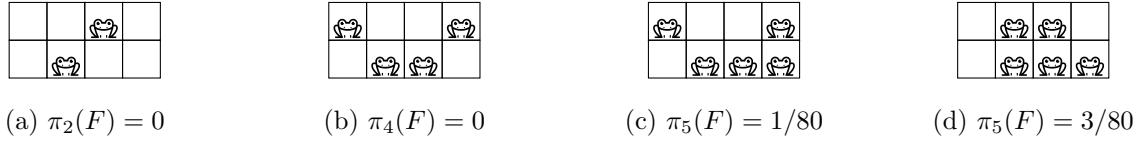
\begin{figure}[ht]
    \begin{center}
        \begin{subfigure}{0.24\textwidth}
            \begin{center}
                \begin{tikzpicture}
                    \clip (-0.5,-0.5) rectangle (2,1.5);
                    \foreach \x in {0,0.5,1,1.5} {
                        \draw ({\x-0.25},-0.25) rectangle (\x+0.25,0.25);
                        \draw ({\x-0.25},0.25) rectangle (\x+0.25,0.75);
                    }
                    \node (f1) at (0.5,0) {$\froggieb$};
                    \node (f2) at (1,0.5) {$\froggieb$};
                \end{tikzpicture}
            \end{center}
            \sfcaption
            \caption{\label{fig:impossible1}$\pi_2(F)=0$}
        \end{subfigure}\hfill
        \begin{subfigure}{0.24\textwidth}
            \begin{center}
                \begin{tikzpicture}
                    \clip (-0.5,-0.5) rectangle (2,1.5);
                    \foreach \x in {0,0.5,1,1.5} {
                        \draw ({\x-0.25},-0.25) rectangle (\x+0.25,0.25);
                        \draw ({\x-0.25},0.25) rectangle (\x+0.25,0.75);
                    }
                    \node (f1) at (0,0.5) {$\froggieb$};
                    \node (f2) at (0.5,0) {$\froggieb$};
                    \node (f3) at (1,0) {$\froggieb$};
                    \node (f4) at (1.5,0.5) {$\froggieb$};
                \end{tikzpicture}
            \end{center}
            \sfcaption
            \caption{\label{fig:impossible2}$\pi_4(F)=0$}
        \end{subfigure}\hfill
        \begin{subfigure}{0.24\textwidth}
            \begin{center}
                \begin{tikzpicture}
                    \clip (-0.5,-0.5) rectangle (2,1.5);
                    \foreach \x in {0,0.5,1,1.5} {
                        \draw ({\x-0.25},-0.25) rectangle (\x+0.25,0.25);
                        \draw ({\x-0.25},0.25) rectangle (\x+0.25,0.75);
                    }
                    \node (f1) at (0,0.5) {$\froggieb$};
                    \node (f2) at (0.5,0) {$\froggieb$};
                    \node (f3) at (1,0) {$\froggieb$};
                    \node (f4) at (1.5,0.5) {$\froggieb$};
                    \node (f5) at (1.5,0) {$\froggieb$};
                \end{tikzpicture}
            \end{center}
            \sfcaption
            \caption{\label{fig:possible1}$\pi_5(F)=1/80$}
        \end{subfigure}
        \begin{subfigure}{0.24\textwidth}
            \begin{center}
                \begin{tikzpicture}
                    \clip (-0.5,-0.5) rectangle (2,1.5);
                    \foreach \x in {0,0.5,1,1.5} {
                        \draw ({\x-0.25},-0.25) rectangle (\x+0.25,0.25);
                        \draw ({\x-0.25},0.25) rectangle (\x+0.25,0.75);
                    }
                    \node (f1) at (0.5,0.5) {$\froggieb$};
                    \node (f2) at (0.5,0) {$\froggieb$};
                    \node (f3) at (1,0) {$\froggieb$};
                    \node (f4) at (1,0.5) {$\froggieb$};
                    \node (f5) at (1.5,0) {$\froggieb$};
                \end{tikzpicture}
            \end{center}
            \sfcaption
            \caption{\label{fig:possible2}$\pi_5(F)=3/80$}
        \end{subfigure}
    \end{center}
    \caption{\label{fig:blind-frog-motivation}Some blind-frog arrangements and their stationary values.}
\end{figure}

Before introducing the frogs with hats, we start with a motivating example.
Consider the blind-frog arrangement $F$ shown in \Cref{fig:impossible1}.
Using $\pi_m$ to denote the stationary distribution of the $m$-blind-frog dynamics, it turns out that $\pi_2(F)=0$... why?
In {\cite[Theorem 20]{bukh_frogs}}, Bukh and Cox argue that any given $m$-blind-frog arrangement is in the support of $\pi_m$ if and only if that arrangement can be reached via a sequence of pokes starting from an arrangement where the $m$ froggies form an interval.
We encourage the reader to check that the arrangement in \Cref{fig:impossible1} is indeed unreachable in this sense.
The arrangement shown in \Cref{fig:impossible2} has $\pi_4(F)=0$ for the same reason, but this is more tedious to verify.
However, the arrangements shown in \Cref{fig:possible1,fig:possible2} both have $\pi_5(F)\neq 0$ and so they must be reachable, even though it may be unclear what sequence of pokes accomplishes this feat.
Despite both being reachable, why is it that the arrangement in \Cref{fig:possible2} is three times more likely to appear than the arrangement in \Cref{fig:possible1}?

All of these questions, and more, will be answered by placing hats on the froggies' heads!

\begin{defn}[Hatted-frog arrangement]
    An $m$-hatted-frog arrangement on the grid $[2]\times[k]$ is a pair of sets $(F,H)$ where $F\in{[2]\times[k]\choose m}$ and $H\subseteq F$ with the following properties:
    \begin{properties}
        \item If $(i,c)\in F$ for some $(i,c)\in[2]\times[k]$, then \emph{exactly one} of $(1,c),(2,c)$ is an element of $H$, and
        \item If $(2,c)\in H$, then $(1,c+1)\notin H$.
    \end{properties}
    We denote the set of all $m$-hatted-frog arrangements on $[2]\times[k]$ by $\hattedFrogs km$.
\end{defn}
To explain the terminology and notation, we imagine that the set $F\in{[2]\times[k]\choose m}=\blindFrogs km$ is the set of positions on the grid occupied by a frog and that $H\subseteq F$ is the set of frogs wearing a hat.
To that end, we refer to each element of $F$ as a `frog' and to each element of $H$ as a `hatted frog'.
We will often use $\froggiehat$ to denote an element of $H$ and $\froggie$ to denote an element of $F\setminus H$.
The whole of $\hattedFrogs 22$ and $\hattedFrogs 23$ can be seen in \Cref{fig:states}.

Informally, an $m$-hatted-frog arrangement is created by placing $m$ frogs (some with hats, some without) in the grid $[2]\times[k]$ so that the following pictures are never seen:
\begin{center}
    \domino{\froggiehat}{\froggiehat}, \domino{\froggieb}{\froggieb}, \domino{\froggieb}{\nofroggie}, \domino{\nofroggie}{\froggieb}, \dominoblock{?}{\froggiehat}{\froggiehat}{?}
\end{center}
where a ``?'' represents a frog that may or may not be present.
In our pictures, the top row corresponds to coordinates of the form $(1,i)$ and the bottom row corresponds to coordinates of the form $(2,i)$.
\medskip

Recalling the rotation map $\rot$ from \Cref{ring_to_grid}, we can naturally apply $\rot$ to each $(F,H)\in\hattedFrogs km$ by applying $\rot$ pointwise to both $F$ and $H$.
It is quick to observe that $\hattedFrogs km$ is invariant under $\rot$.
\medskip

At this point, the stationary values in \Cref{fig:blind-frog-motivation} begin to make sense.
There are no valid ways to place hats on the frogs in \Cref{fig:impossible1,fig:impossible2}, whereas there is a single way to place hats on the frogs in \Cref{fig:possible1} and three ways to place hats on the frogs in \Cref{fig:possible2} (it turns out that $\abs{\hattedFrogs 45}=80$).

There is an obvious projection $\doff\colon\hattedFrogs km\to\blindFrogs km$ which simply removes the hats from all frogs; that is $\doff(F,H)=F$.
In order to explain the numbers in \Cref{fig:blind-frog-motivation}, our goal is to create a Markov chain on $\hattedFrogs km$ such that $\doff$ realizes a coupling to the $m$-blind-frog dynamics.
If all goes well, this Markov chain will wind up having a uniform stationary distribution.

\subsection{The hatted-frog process and dynamics}\label{hat_process}

In this section, we build a Markov chain on $\hattedFrogs km$ which projects to the $m$-blind-frog dynamics for $W=12\cdots kk\cdots 21$.
Given a hatted-frog arrangement $\hat F\in\hattedFrogs km$, the \emph{$m$-hatted-frog process} evolves as follows:
\begin{enumerate}
    \item\label{hatted_defn} The monster which lives in the pond selects some letter $a\in\bet$.
        If $a\notin[k]$, then nothing happens; otherwise, the monster pokes column $a$.
        This agitates the frogs in column $a$ and makes them want to jump away.
    \item An agitated frog will hop to the very next square in the clockwise order.
        If there were initially two agitated frogs, then the one which was \emph{not} originally wearing a hat hops first.
    \item When an agitated frog lands on a new square, it calms down and immediately puts on a hat, potentially stealing this hat from the current occupant or from the occupant of the other square in this column, should either occupant exist.
        If the square was previously occupied, then that frog becomes agitated and will hop \emph{immediately} on the next step.
        In particular, this newly agitated frog hops before any previously agitated frogs.

        This repeats until all frogs are content once more.
\end{enumerate}

For an $m$-hatted-frog arrangement $\hat F\in\hattedFrogs km$ and a letter $a\in\bet$, let $\hat Fa$ be the hatted-frog arrangement that results from starting with $\hat F$, poking column $a$ (if $a\in[k]$) and waiting for the ensuing frenzy to settle.
Additionally, define $\hop(\hat F,a)$ to be the number of frogs which hopped in the transition from $\hat F$ to $\hat Fa$.
Note that if $a\notin[k]$ or $\hat F$ has no frogs in column $a$, then $\hop(\hat F,a)=0$.
\medskip

Examples of the hatted-frog process can be seen in \Cref{fig:hattedProcess}.
\begin{figure}
    \begin{subfigure}{\textwidth}
        \begin{center}
            \begin{tikzpicture}
                \clip (-0.5,-1) rectangle (2,1.5);
                \foreach \x in {0,0.5,1,1.5} {
                    \draw ({\x-0.25},-0.25) rectangle (\x+0.25,0.25);
                    \draw ({\x-0.25},0.25) rectangle (\x+0.25,0.75);
                }

                \node (f1) at (0,0.5) {$\froggiehat$};
                \node (f2) at (0.5,0.5) {$\froggiehat$};
                \node (f3) at (1.5,0.5) {$\froggieb$};
                \node (f4) at (0,0) {$\froggieb$};
                \node (f3) at (1,0) {$\froggiehat$};
                \node (f4) at (1.5,0) {$\froggiehat$};
            \end{tikzpicture}
            \begin{tikzpicture}
                \clip (-0.5,-1) rectangle (2,1.5);
                \foreach \x in {0,0.5,1,1.5} {
                    \draw ({\x-0.25},-0.25) rectangle (\x+0.25,0.25);
                    \draw ({\x-0.25},0.25) rectangle (\x+0.25,0.75);
                }

                \node (f1) at (0,0.5) {$\froggiehat$};
                \node (f2) at (0.5,1) {$\froggiehat$};
                \node (f3) at (1.5,0.5) {$\froggieb$};
                \node (f4) at (0,0) {$\froggieb$};
                \node (f3) at (1,0) {$\froggiehat$};
                \node (f4) at (1.5,0) {$\froggiehat$};
            \end{tikzpicture}
            \begin{tikzpicture}
                \clip (-0.5,-1) rectangle (2,1.5);
                \foreach \x in {0,0.5,1,1.5} {
                    \draw ({\x-0.25},-0.25) rectangle (\x+0.25,0.25);
                    \draw ({\x-0.25},0.25) rectangle (\x+0.25,0.75);
                }

                \node (f1) at (0,0.5) {$\froggiehat$};
                \node (f2) at (1,0.5) {$\froggiehat$};
                \node (f3) at (1.5,0.5) {$\froggieb$};
                \node (f4) at (0,0) {$\froggieb$};
                \node (f3) at (1,0) {$\froggieb$};
                \node (f4) at (1.5,0) {$\froggiehat$};
            \end{tikzpicture}
        \end{center}
        \sfcaption
        \caption{\label{fig:hattedProcess:a}
            The monster pokes column $2$!
            The left-most image is $\hat F$, the right-most image is $\hat F2$ and $\hop(\hat F,2)=1$.
        }
    \end{subfigure}
    \sfbetween

    \begin{subfigure}{\textwidth}
        \begin{center}
            \begin{tikzpicture}
                \clip (-0.5,-1) rectangle (2,1.5);
                \foreach \x in {0,0.5,1,1.5} {
                    \draw ({\x-0.25},-0.25) rectangle (\x+0.25,0.25);
                    \draw ({\x-0.25},0.25) rectangle (\x+0.25,0.75);
                }

                \node (f1) at (0.5,0.5) {$\froggiehat$};
                \node (f2) at (1.5,0.5) {$\froggieb$};
                \node (f3) at (1,0) {$\froggiehat$};
                \node (f4) at (1.5,0) {$\froggiehat$};
            \end{tikzpicture}
            \begin{tikzpicture}
                \clip (-0.5,-1) rectangle (2,1.5);
                \foreach \x in {0,0.5,1,1.5} {
                    \draw ({\x-0.25},-0.25) rectangle (\x+0.25,0.25);
                    \draw ({\x-0.25},0.25) rectangle (\x+0.25,0.75);
                }

                \node (f1) at (0.5,0.5) {$\froggiehat$};
                \node (f2) at (1.5,1) {$\froggie$};
                \node (f3) at (1,0) {$\froggiehat$};
                \node (f4) at (1.5,-0.5) {$\froggiehat$};
            \end{tikzpicture}
            \begin{tikzpicture}
                \clip (-0.5,-1) rectangle (2,1.5);
                \foreach \x in {0,0.5,1,1.5} {
                    \draw ({\x-0.25},-0.25) rectangle (\x+0.25,0.25);
                    \draw ({\x-0.25},0.25) rectangle (\x+0.25,0.75);
                }

                \node (f1) at (0.5,0.5) {$\froggiehat$};
                \node (f2) at (1.5,0) {$\froggiehat$};
                \node (f3) at (1,0) {$\froggiehat$};
                \node (f4) at (1.5,-0.5) {$\froggiehat$};
            \end{tikzpicture}
            \begin{tikzpicture}
                \clip (-0.5,-1) rectangle (2,1.5);
                \foreach \x in {0,0.5,1,1.5} {
                    \draw ({\x-0.25},-0.25) rectangle (\x+0.25,0.25);
                    \draw ({\x-0.25},0.25) rectangle (\x+0.25,0.75);
                }

                \node (f1) at (0.5,0.5) {$\froggiehat$};
                \node (f2) at (1.5,0) {$\froggiehat$};
                \node (f3) at (1,-0.5) {$\froggie$};
                \node (f4) at (1,0) {$\froggiehat$};
            \end{tikzpicture}
            \begin{tikzpicture}
                \clip (-0.5,-1) rectangle (2,1.5);
                \foreach \x in {0,0.5,1,1.5} {
                    \draw ({\x-0.25},-0.25) rectangle (\x+0.25,0.25);
                    \draw ({\x-0.25},0.25) rectangle (\x+0.25,0.75);
                }

                \node (f1) at (0.5,0.5) {$\froggieb$};
                \node (f2) at (1.5,0) {$\froggiehat$};
                \node (f3) at (0.5,0) {$\froggiehat$};
                \node (f4) at (1,0) {$\froggiehat$};
            \end{tikzpicture}
        \end{center}
        \sfcaption
        \caption{\label{fig:hattedProcess:b}
            The monster pokes column $4$!
            The left-most image is $\hat F$, the right-most image is $\hat F4$ and $\hop(\hat F,4)=3$.
        }
    \end{subfigure}
    \sfbetween

    \begin{subfigure}{\textwidth}
        \begin{center}
            \begin{tikzpicture}
                \clip (-0.5,-1) rectangle (2,1.5);
                \foreach \x in {0,0.5,1,1.5} {
                    \draw ({\x-0.25},-0.25) rectangle (\x+0.25,0.25);
                    \draw ({\x-0.25},0.25) rectangle (\x+0.25,0.75);
                }

                \node (f1) at (0,0.5) {$\froggiehat$};
                \node (f2) at (0.5,0.5) {$\froggieb$};
                \node (f3) at (1,0.5) {$\froggieb$};
                \node (f4) at (0,0) {$\froggieb$};
                \node (f5) at (0.5,0) {$\froggiehat$};
                \node (f6) at (1,0) {$\froggiehat$};
                \node (f7) at (1.5,0) {$\froggiehat$};
            \end{tikzpicture}
            \begin{tikzpicture}
                \clip (-0.5,-1) rectangle (2,1.5);
                \foreach \x in {0,0.5,1,1.5} {
                    \draw ({\x-0.25},-0.25) rectangle (\x+0.25,0.25);
                    \draw ({\x-0.25},0.25) rectangle (\x+0.25,0.75);
                }

                \node (f1) at (0,0.5) {$\froggiehat$};
                \node (f2) at (0.5,1) {$\froggie$};
                \node (f3) at (1,0.5) {$\froggieb$};
                \node (f4) at (0,0) {$\froggieb$};
                \node (f5) at (0.5,-0.5) {$\froggiehat$};
                \node (f6) at (1,0) {$\froggiehat$};
                \node (f7) at (1.5,0) {$\froggiehat$};
            \end{tikzpicture}
            \begin{tikzpicture}
                \clip (-0.5,-1) rectangle (2,1.5);
                \foreach \x in {0,0.5,1,1.5} {
                    \draw ({\x-0.25},-0.25) rectangle (\x+0.25,0.25);
                    \draw ({\x-0.25},0.25) rectangle (\x+0.25,0.75);
                }

                \node (f1) at (0,0.5) {$\froggiehat$};
                \node (f2) at (1,0.5) {$\froggiehat$};
                \node (f3) at (1,1) {$\froggie$};
                \node (f4) at (0,0) {$\froggieb$};
                \node (f5) at (0.5,-0.5) {$\froggiehat$};
                \node (f6) at (1,0) {$\froggieb$};
                \node (f7) at (1.5,0) {$\froggiehat$};
            \end{tikzpicture}
            \begin{tikzpicture}
                \clip (-0.5,-1) rectangle (2,1.5);
                \foreach \x in {0,0.5,1,1.5} {
                    \draw ({\x-0.25},-0.25) rectangle (\x+0.25,0.25);
                    \draw ({\x-0.25},0.25) rectangle (\x+0.25,0.75);
                }

                \node (f1) at (0,0.5) {$\froggiehat$};
                \node (f2) at (1,0.5) {$\froggiehat$};
                \node (f3) at (1.5,0.5) {$\froggiehat$};
                \node (f4) at (0,0) {$\froggieb$};
                \node (f5) at (0.5,-0.5) {$\froggiehat$};
                \node (f6) at (1,0) {$\froggieb$};
                \node (f7) at (1.5,0) {$\froggieb$};
            \end{tikzpicture}
            \begin{tikzpicture}
                \clip (-0.5,-1) rectangle (2,1.5);
                \foreach \x in {0,0.5,1,1.5} {
                    \draw ({\x-0.25},-0.25) rectangle (\x+0.25,0.25);
                    \draw ({\x-0.25},0.25) rectangle (\x+0.25,0.75);
                }

                \node (f1) at (0,0.5) {$\froggieb$};
                \node (f2) at (1,0.5) {$\froggiehat$};
                \node (f3) at (1.5,0.5) {$\froggiehat$};
                \node (f4) at (0,-0.5) {$\froggie$};
                \node (f5) at (0,0) {$\froggiehat$};
                \node (f6) at (1,0) {$\froggieb$};
                \node (f7) at (1.5,0) {$\froggieb$};
            \end{tikzpicture}\\
            \begin{tikzpicture}
                \clip (-0.5,-1) rectangle (2,1.5);
                \foreach \x in {0,0.5,1,1.5} {
                    \draw ({\x-0.25},-0.25) rectangle (\x+0.25,0.25);
                    \draw ({\x-0.25},0.25) rectangle (\x+0.25,0.75);
                }

                \node (f1) at (0,1) {$\froggie$};
                \node (f2) at (1,0.5) {$\froggiehat$};
                \node (f3) at (1.5,0.5) {$\froggiehat$};
                \node (f4) at (0,0.5) {$\froggiehat$};
                \node (f5) at (0,0) {$\froggieb$};
                \node (f6) at (1,0) {$\froggieb$};
                \node (f7) at (1.5,0) {$\froggieb$};
            \end{tikzpicture}
            \begin{tikzpicture}
                \clip (-0.5,-1) rectangle (2,1.5);
                \foreach \x in {0,0.5,1,1.5} {
                    \draw ({\x-0.25},-0.25) rectangle (\x+0.25,0.25);
                    \draw ({\x-0.25},0.25) rectangle (\x+0.25,0.75);
                }

                \node (f1) at (0.5,0.5) {$\froggiehat$};
                \node (f2) at (1,0.5) {$\froggiehat$};
                \node (f3) at (1.5,0.5) {$\froggiehat$};
                \node (f4) at (0,0.5) {$\froggiehat$};
                \node (f5) at (0,0) {$\froggieb$};
                \node (f6) at (1,0) {$\froggieb$};
                \node (f7) at (1.5,0) {$\froggieb$};
            \end{tikzpicture}
        \end{center}
        \sfcaption
        \caption{\label{fig:hattedProcess:c}
            The monster pokes column $2$!
            The top-left image is $\hat F$, the bottom-right image is $\hat F2$ and $\hop(\hat F,2)=5$.
        }
    \end{subfigure}
    \sfbetween

    \begin{subfigure}{\textwidth}
        \begin{center}
            \begin{tikzpicture}
                \clip (-0.5,-1) rectangle (2,1.5);
                \foreach \x in {0,0.5,1,1.5} {
                    \draw ({\x-0.25},-0.25) rectangle (\x+0.25,0.25);
                    \draw ({\x-0.25},0.25) rectangle (\x+0.25,0.75);
                }

                \node (f1) at (0,0.5) {$\froggiehat$};
                \node (f2) at (0.5,0.5) {$\froggiehat$};
                \node (f3) at (1,0.5) {$\froggieb$};
                \node (f4) at (0,0) {$\froggieb$};
                \node (f5) at (0.5,0) {$\froggieb$};
                \node (f6) at (1,0) {$\froggiehat$};
                \node (f7) at (1.5,0) {$\froggiehat$};
            \end{tikzpicture}
            \begin{tikzpicture}
                \clip (-0.5,-1) rectangle (2,1.5);
                \foreach \x in {0,0.5,1,1.5} {
                    \draw ({\x-0.25},-0.25) rectangle (\x+0.25,0.25);
                    \draw ({\x-0.25},0.25) rectangle (\x+0.25,0.75);
                }

                \node (f1) at (0,0.5) {$\froggiehat$};
                \node (f2) at (0.5,0.5) {$\froggiehat$};
                \node (f3) at (1,1) {$\froggie$};
                \node (f4) at (0,0) {$\froggieb$};
                \node (f5) at (0.5,0) {$\froggieb$};
                \node (f6) at (1,-0.5) {$\froggiehat$};
                \node (f7) at (1.5,0) {$\froggiehat$};
            \end{tikzpicture}
            \begin{tikzpicture}
                \clip (-0.5,-1) rectangle (2,1.5);
                \foreach \x in {0,0.5,1,1.5} {
                    \draw ({\x-0.25},-0.25) rectangle (\x+0.25,0.25);
                    \draw ({\x-0.25},0.25) rectangle (\x+0.25,0.75);
                }

                \node (f1) at (0,0.5) {$\froggiehat$};
                \node (f2) at (0.5,0.5) {$\froggiehat$};
                \node (f3) at (1.5,0.5) {$\froggiehat$};
                \node (f4) at (0,0) {$\froggieb$};
                \node (f5) at (0.5,0) {$\froggieb$};
                \node (f6) at (1,-0.5) {$\froggiehat$};
                \node (f7) at (1.5,0) {$\froggieb$};
            \end{tikzpicture}
            \begin{tikzpicture}
                \clip (-0.5,-1) rectangle (2,1.5);
                \foreach \x in {0,0.5,1,1.5} {
                    \draw ({\x-0.25},-0.25) rectangle (\x+0.25,0.25);
                    \draw ({\x-0.25},0.25) rectangle (\x+0.25,0.75);
                }

                \node (f1) at (0,0.5) {$\froggiehat$};
                \node (f2) at (0.5,0.5) {$\froggieb$};
                \node (f3) at (1.5,0.5) {$\froggiehat$};
                \node (f4) at (0,0) {$\froggieb$};
                \node (f5) at (0.5,-0.5) {$\froggie$};
                \node (f6) at (0.5,0) {$\froggiehat$};
                \node (f7) at (1.5,0) {$\froggieb$};
            \end{tikzpicture}
            \begin{tikzpicture}
                \clip (-0.5,-1) rectangle (2,1.5);
                \foreach \x in {0,0.5,1,1.5} {
                    \draw ({\x-0.25},-0.25) rectangle (\x+0.25,0.25);
                    \draw ({\x-0.25},0.25) rectangle (\x+0.25,0.75);
                }

                \node (f1) at (0,0.5) {$\froggieb$};
                \node (f2) at (0.5,0.5) {$\froggieb$};
                \node (f3) at (1.5,0.5) {$\froggiehat$};
                \node (f4) at (0,-0.5) {$\froggie$};
                \node (f5) at (0,0) {$\froggiehat$};
                \node (f6) at (0.5,0) {$\froggiehat$};
                \node (f7) at (1.5,0) {$\froggieb$};
            \end{tikzpicture}\\
            \begin{tikzpicture}
                \clip (-0.5,-1) rectangle (2,1.5);
                \foreach \x in {0,0.5,1,1.5} {
                    \draw ({\x-0.25},-0.25) rectangle (\x+0.25,0.25);
                    \draw ({\x-0.25},0.25) rectangle (\x+0.25,0.75);
                }

                \node (f1) at (0,1) {$\froggie$};
                \node (f2) at (0.5,0.5) {$\froggieb$};
                \node (f3) at (1.5,0.5) {$\froggiehat$};
                \node (f4) at (0,0.5) {$\froggiehat$};
                \node (f5) at (0,0) {$\froggieb$};
                \node (f6) at (0.5,0) {$\froggiehat$};
                \node (f7) at (1.5,0) {$\froggieb$};
            \end{tikzpicture}
            \begin{tikzpicture}
                \clip (-0.5,-1) rectangle (2,1.5);
                \foreach \x in {0,0.5,1,1.5} {
                    \draw ({\x-0.25},-0.25) rectangle (\x+0.25,0.25);
                    \draw ({\x-0.25},0.25) rectangle (\x+0.25,0.75);
                }

                \node (f1) at (0.5,0.5) {$\froggiehat$};
                \node (f2) at (0.5,1) {$\froggie$};
                \node (f3) at (1.5,0.5) {$\froggiehat$};
                \node (f4) at (0,0.5) {$\froggiehat$};
                \node (f5) at (0,0) {$\froggieb$};
                \node (f6) at (0.5,0) {$\froggieb$};
                \node (f7) at (1.5,0) {$\froggieb$};
            \end{tikzpicture}
            \begin{tikzpicture}
                \clip (-0.5,-1) rectangle (2,1.5);
                \foreach \x in {0,0.5,1,1.5} {
                    \draw ({\x-0.25},-0.25) rectangle (\x+0.25,0.25);
                    \draw ({\x-0.25},0.25) rectangle (\x+0.25,0.75);
                }

                \node (f1) at (0.5,0.5) {$\froggiehat$};
                \node (f2) at (1,0.5) {$\froggiehat$};
                \node (f3) at (1.5,0.5) {$\froggiehat$};
                \node (f4) at (0,0.5) {$\froggiehat$};
                \node (f5) at (0,0) {$\froggieb$};
                \node (f6) at (0.5,0) {$\froggieb$};
                \node (f7) at (1.5,0) {$\froggieb$};
            \end{tikzpicture}
        \end{center}
        \sfcaption
        \caption{\label{fig:hattedProcess:d}
            The monster pokes column $3$!
            The top-left image is $\hat F$, the bottom-right is $\hat F3$ and $\hop(\hat F,3)=6$.
        }
    \end{subfigure}
    \caption{\label{fig:hattedProcess}
        Examples of the \hyperref[hatted_defn]{hatted-frog process}.
        An agitated frog is drawn just above/below the square that it occupies.
        We emphasize the hatted poked frog by continuing to draw a hat on its head until it first hops.
    }
\end{figure}

\medskip

In the hatted-frog process, if two frogs are initially agitated, then the one that was originally wearing a hat waits until it is the only agitated frog before it hops.
The astute reader may realize that this requirement is somewhat arbitrary.
Indeed, this choice makes a difference only when there are exactly $2k$ frogs (in which case it is easy to verify that $\hop(F,a)=2k$ for all $a\in[k]$).
Nevertheless, enforcing this hopping order will aid our analysis of the hatted-frog process.
\medskip

Recalling the rotation map $\rot$ from \Cref{ring_to_grid}, it is clear that:

\begin{observation}\label{rotatehats}
    Fix a letter $c\in[k]$ and set $c'=k+1-c$.
    Then for any $\hat F\in\hattedFrogs km$, we have $\rot(\hat F)\in\hattedFrogs km$ and $\rot(\hat Fc)=\rot(\hat F)c'$.
\end{observation}

We next show that the hatted-frog process extends the blind-frog process.
Recall that $\doff$ is the projection which simply removes the frogs' hats.

\begin{theorem}\label{tohats}
    For any letter $a\in\bet$ and any hatted-frog arrangement $\hat F\in\hattedFrogs km$,
    \begin{align*}
        \doff(\hat Fa) &=\doff(\hat F)a,\qquad\text{and}\\
        \hop(\hat F,a) &= \hop(\doff(\hat F),a).
    \end{align*}
\end{theorem}

In order to prove this, we introduce a pair of definitions and a lemma which will be used throughout the next section as well.

\begin{defn}
    Fix a square $\square\in[2]\times[k]$.
    \begin{itemize}
        \item The column containing $\square$ is denoted by $\column(\square)$; that is to say, $\column(r,c)\eqdef [2]\times\{c\}$.
        \item The square opposite to $\square$ is denoted by $\opp(\square)$; that is to say, $\opp(r,c)\eqdef (3-r,c)$.
        \item $\square^+$ denotes the square directly succeeding $\square$ in the clockwise order and $\square^-$ denotes the square directly preceding $\square$ in the clockwise order.
            For example, $(1,1)^+=(1,2)$ and $(1,1)^-=(2,1)$ and $(2,k)^+=(2,k-1)$ and $(2,k)^-=(1,k)$.
    \end{itemize}
\end{defn}
Note that $\column$ and $\opp$ both commute with $\rot$ and that $\rot(\square^{\pm})=\rot(\square)^{\pm}$.

Recall that we use $\Omega+x$ as shorthand for $\Omega\cup\{x\}$ and use use $\Omega-x$ as shorthand for $\Omega\setminus\{x\}$.

\begin{lemma}\label{nextStep}
    Suppose that $(F,H)$ is a hatted-frog arrangement on $[2]\times[k]$ and fix any $\square\in[2]\times[k]$.
    If either $\square^-\in H$ or $F\cap\column(\square^-)=\varnothing$, then $\bigl(F+\square,\ H-\opp(\square)+\square\bigr)$ is also a hatted-frog arrangement.
\end{lemma}
\begin{proof}
    Using rotational symmetry, we may suppose that $\square=(1,c)$ for some $c\in[k]$.
    It is quick to observe that the only way for $\bigl(F+\square,H-\opp(\square)+\square\bigr)$ to fail to be a hatted-frog arrangement is if $c\geq 2$ and $(2,c-1)\in H-\opp(\square)+\square$.
    However, by assumption, either $\square^-=(1,c-1)\in H$ or $H$ is disjoint from $\column(\square^-)=[2]\times\{c-1\}$.
    In either case $(2,c-1)\notin H-\opp(\square)+\square$.
\end{proof}

\begin{proof}[Proof of \Cref{tohats}]
    Observe that the only heed given to the presence of the hats in the hatted-frog process is to determine which of the agitated frogs has priority in moving.
    However, in the blind-frog process, the order in which the froggies hop is irrelevant.
    Therefore, we will have established the theorem provided that we show that $\hat Fa\in\hattedFrogs km$ for each $\hat F\in\hattedFrogs km$ and each $a\in\bet$, i.e.\ we need only show that the hatted-frog process really is a process on $\hattedFrogs km$.
    Of course, if $a\notin[k]$, then $\hat F a=\hat F$, so we may suppose that $a\in[k]$.
    \medskip

    Set $h=\hop(\hat F,a)$.
    Fix $i\in[h]$ and consider the state of the hatted frogs after $i$ hops in the transition from $\hat F$ to $\hat F a$ have taken place.
    We can encode this state in a triple $(F_i,H_i,A_i)$ where $A_i\in{[2]\times[k]\choose\leq 2}$ is the set of currently agitated frogs, $F_i\in{[2]\times[k]\choose \leq m}$ is the set of un-agitated frogs and $H_i\subseteq F_i$ is the set of those un-agitated frogs wearing a hat.
    Note that $\abs{F_i}+\abs{A_i}=m$ and that $F_i$ and $A_i$ may intersect.
    We additionally define $(F_0,H_0,A_0)$ to be the state immediately after column $a$ is poked but before any frogs have hopped.
    In particular, if $\hat F=(F,H)$, then $A_0=F\cap\bigl([2]\times\{a\}\bigr)$, $F_0=F\setminus A_0$ and $H_0=H\setminus A_0$.
    Certainly, $(F_0,H_0)$ is a hatted-frog arrangement.

    By definition, $A_h=\varnothing$ and $\hat F a=(F_h,H_h)$.
    Therefore, the goal is to show that $(F_h,H_h)$ is indeed a hatted-frog arrangement.

    Suppose that $\froggie\in[2]\times[k]$ is the $i$th frog to hop.
    Then $\froggie\in A_{i-1}$ and
    \begin{itemize}
        \item $F_i=F_{i-1}+\froggie^+$, and
        \item $H_i= H_{i-1}-\opp(\froggie^+)+\froggie^+$, and
        \item If $\froggie^+\notin F$, then $A_i=A_{i-1}-\froggie$, and
        \item If $\froggie^+\in F$, then $A_i=A_{i-1}-\froggie+\froggie^+$.
    \end{itemize}
    Using this information, it is not difficult to see that for each $i\in[0,h]$ and each $\froggie\in A_i$, either $\froggie\in H_i$ or $F_i\cap\column(\froggie)=\varnothing$.
    Therefore, a routine induction using \Cref{nextStep} shows that $(F_i,H_i)$ is indeed a hatted-frog arrangement for each $i\in[0,h]$, which concludes the proof.
\end{proof}

Naturally, we can build a Markov chain from the hatted-frog process just as we have done from the frog and blind-frog processes.
Indeed, fix integers $k,m$ with $0\leq m\leq 2k$ and an alphabet $\bet$.
Starting with a hatted-frog arrangement $\hat F_0\in\hattedFrogs km$ set $\hat F_n=\hat F_0r_1r_2\cdots r_n$ where the letters $r_1,r_2,\ldots$ are chosen independently and uniformly from $\bet$.
Then the sequence $\hat F_0,\hat F_1,\hat F_2,\ldots$ is a Markov chain, which we call the \emph{$(k,\bet)$ $m$-hatted-frog dynamics}.

\Cref{tohats} tells us that we can couple the $(12\cdots kk\cdots 21,\bet)$ $m$-blind-frog dynamics and the $(k,\bet)$ $m$-hatted-frog dynamics by simply poking the same letter in each chain at each time step.

\begin{corollary}\label{coupling}
    Fix $W=12\cdots kk\cdots 21$, some alphabet $\bet\supseteq[k]$ and an integer $0\leq m\leq 2k$.
    If $\pi$ denotes the stationary distribution of the $(W,\bet)$ $m$-blind-frog dynamics and $\hat\pi$ is \emph{any} stationary distribution of the $(k,\bet)$ $m$-hatted-frog dynamics, then
    \begin{enumerate}
        \item $\hat F\sim\hat\pi\implies \doff(\hat F)\sim\pi$, and
        \item $\displaystyle \E_{F\sim\pi} \hop(F,a) = \E_{\scalebox{0.65}{$\hat F\!\sim\!\hat\pi$}} \hop(\hat F,a)$ for each $a\in\bet$.
    \end{enumerate}
\end{corollary}

With this coupling in hand, the next step is to locate a stationary distribution of the hatted-frog dynamics.

\begin{theorem}\label{uniform}
    If $\bet\supseteq[k]$, then the $(k,\bet)$ $m$-hatted-frog dynamics admits a uniform stationary distribution.
\end{theorem}
Note that we prove only that the uniform distribution is \emph{a} stationary distribution, not that it is unique.\footnote{It should be possible to prove that the stationary distribution is indeed unique, but it is not necessary for our arguments.}
Of course, \Cref{coupling} and \Cref{uniform} together imply that the stationary distribution of the $(12\cdots kk\cdots 21,\bet)$ $m$-blind-frog dyanimics where $\bet\supseteq[k]$, is given by
\[
    \pi(F)={\abs{\doff^{-1}(F)}\over\abs{\hattedFrogs km}}.
\]
The whole of \Cref{sec:stationary} is dedicated to the proof of \Cref{uniform}.

\subsection{The hatted-frog dynamics admits a uniform stationary distribution}\label{sec:stationary}
To prove \Cref{uniform}, it suffices to show that the state-graph of the $m$-hatted-frog process is regular.

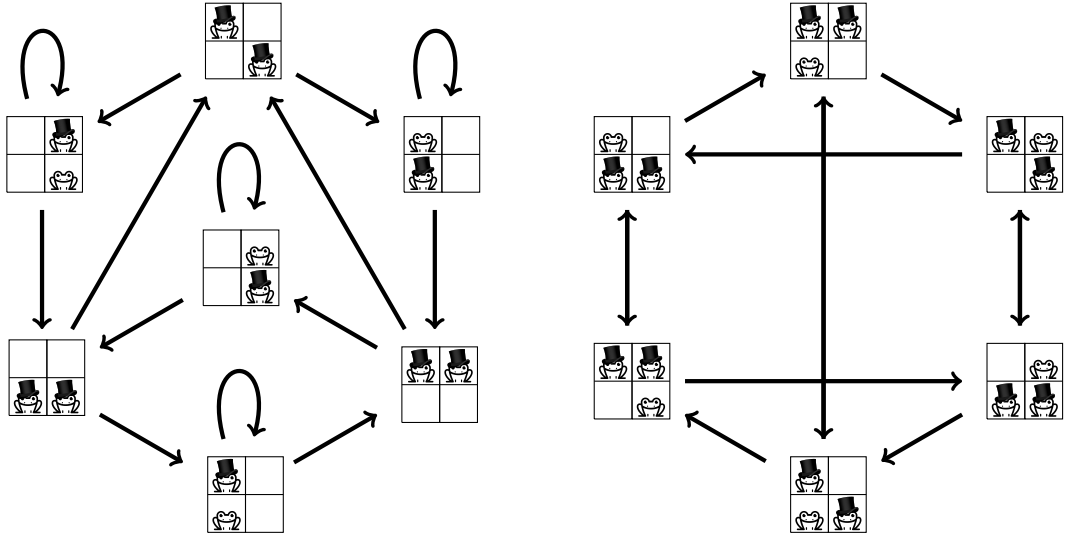
\begin{figure}
    \begin{center}
        \def\s{3}
        \begin{tikzpicture}
            \node (f0) at (-30:\s) {\dominoblock{\froggiehat}{}{\froggiehat}{}};
            \node (f1) at (210:\s) {\dominoblock{}{\froggiehat}{}{\froggiehat}};
            \node (f2) at (-90:\s) {\dominoblock{\froggiehat}{\froggieb}{}{}};
            \node (f3) at (30:\s) {\dominoblock{\froggieb}{\froggiehat}{}{}};
            \node (f4) at (150:\s) {\dominoblock{}{}{\froggiehat}{\froggieb}};
            \node (f5) at (-90:0) {\dominoblock{}{}{\froggieb}{\froggiehat}};
            \node (f6) at (90:\s) {\dominoblock{\froggiehat}{}{}{\froggiehat}};

            \draw[ultra thick, ->] (f0) -- (f5);
            \draw[ultra thick, ->] (f0) -- (f6);
            \draw[ultra thick, ->] (f1) -- (f6);
            \draw[ultra thick, ->] (f1) -- (f2);
            \draw[ultra thick, ->] (f2) -- (f0);
            \draw[ultra thick, ->] (f2) edge[loop above] (f2);
            \draw[ultra thick, ->] (f3) -- (f0);
            \draw[ultra thick, ->] (f3) edge[loop above] (f3);
            \draw[ultra thick, ->] (f4) edge[loop above] (f4);
            \draw[ultra thick, ->] (f4) -- (f1);
            \draw[ultra thick, ->] (f5) edge[loop above] (f5);
            \draw[ultra thick, ->] (f5) -- (f1);
            \draw[ultra thick, ->] (f6) -- (f4);
            \draw[ultra thick, ->] (f6) -- (f3);
        \end{tikzpicture}\hfil
        \begin{tikzpicture}
            \node (f0) at (90:\s) {\dominoblock{\froggiehat}{\froggieb}{\froggiehat}{}};
            \node (f1) at (210:\s) {\dominoblock{\froggiehat}{}{\froggiehat}{\froggieb}};
            \node (f2) at (30:\s) {\dominoblock{\froggiehat}{}{\froggieb}{\froggiehat}};
            \node (f3) at (-90:\s) {\dominoblock{\froggiehat}{\froggieb}{}{\froggiehat}};
            \node (f4) at (150:\s) {\dominoblock{\froggieb}{\froggiehat}{}{\froggiehat}};
            \node (f5) at (-30:\s) {\dominoblock{}{\froggiehat}{\froggieb}{\froggiehat}};

            \draw[ultra thick, ->] (f0) -- (f2);
            \draw[ultra thick, ->] (f0) -- (f3);
            \draw[ultra thick, ->] (f1) -- (f5);
            \draw[ultra thick, ->] (f1) -- (f4);
            \draw[ultra thick, ->] (f2) -- (f5);
            \draw[ultra thick, ->] (f2) -- (f4);
            \draw[ultra thick, ->] (f3) -- (f1);
            \draw[ultra thick, ->] (f3) -- (f0);
            \draw[ultra thick, ->] (f4) -- (f1);
            \draw[ultra thick, ->] (f4) -- (f0);
            \draw[ultra thick, ->] (f5) -- (f2);
            \draw[ultra thick, ->] (f5) -- (f3);
        \end{tikzpicture}
    \end{center}
    \caption{\label{fig:states}
        Two examples of the state-graph of the hatted-frog dynamics on the $2\times 2$ grid using the alphabet $\bet=[2]$.
        Observe that each graph is $2$-regular.
    }
\end{figure}

Let $G$ be the directed graph whose vertex set is $\hattedFrogs km$ where $\hat E\to \hat F$ in $G$ if there is some $a\in\bet$ for which $\hat F=\hat Ea$.
An example of such a graph is shown in \Cref{fig:states}.
Note that $G$ may have multi-edges (if $\hat F=\hat E a=\hat E b$ for some $a\neq b\in\bet$) and loops (if $\hat F=\hat Fa$ for some $a\in\bet$).\footnote{It turns out that multi-edges (except for multi-loops) are not possible and that loops occur only when an empty column is poked. In fact, these observations hold (with one exception) for any instance of the blind-frog dynamics. However, our arguments do not benefit from such knowledge and so we do not spend the time to prove this. The interested reader may find it enlightening to prove such a fact, though.}
With these considerations, $\deg^+(\hat F)=\abs\bet$ for each $\hat F\in\hattedFrogs km$.
By definition, the $m$-hatted-frog dynamics is precisely a random walk on $G$ where each edge is traversed with equal probability.

\begin{theorem}\label{regular}
    If $\bet\supseteq[k]$, then $\deg^+(\hat F)=\deg^-(\hat F)=\abs\bet$ for every $\hat F\in\hattedFrogs km$.
\end{theorem}
We already know that $\deg^+(\hat F)=\abs\bet$ from the earlier observation and so our efforts will be focused on proving that $\deg^-(\hat F)=\abs\bet$.

Note that if the state graph is regular, then the uniform distribution is a stationary distribution.
This standard fact follows immediately from the stationary equation.
Therefore, \Cref{regular} implies \Cref{uniform}.

Establishing \Cref{regular} will comprise the remainder of this section and will require a careful analysis of the intermediate steps within each transition in the hatted-frog process.
\medskip

Before diving into the details, we outline the motivating idea.
Consider a hatted-frog arrangement $\hat F=(F,H)\in\hattedFrogs km$.
Observe that the number of non-empty columns in $\hat F$ is precisely $\abs H$.
Therefore, for a letter $a\in\bet$, if either $a\notin[k]$ or there are no frogs in column $a$, then $\hat Fa=\hat F$.
This accounts for $\abs\bet-\abs H$ many in-edges to $\hat F$.
As such, we need to locate exactly $\abs H$ many additional in-edges; each of these in-edges will have the property that some frogs hopped in the transition.
With this set-up, it is natural to try to find a bijection between $H$ and pairs $(\hat E,a)\in\hattedFrogs km\times[k]$ for which $\hat Ea=\hat F$ and $\hat E$ contains some frog in column $a$.
Of course, since $\hat E$ has some frog in column $a$, it has a hatted-frog in column $a$.
Our goal is to track this special hatted-frog as it hops along in such a way that the ending-position of the special hatted-frog in $\hat F$ uniquely identifies the pair $(\hat E,a)$.

\subsubsection{Preliminaries}
In this section, we lay out a number of definitions and results that will be crucial to establishing that the hatted-frog dynamics admits a uniform stationary distribution.

To begin, recall that if $(F,H)$ is a hatted-frog arrangement, then $H$ cannot contain both $(2,c)$ and $(1,c+1)$ for any $c$.
We note the following rephrasing of this property which will streamline some arguments.
\begin{observation}\label{oppPlus}
    If $(F,H)\in\hattedFrogs km$ and $\froggiehat\in H$, then either $\opp(\froggiehat)^+=\froggiehat$ or $\opp(\froggiehat)^+\notin H$.
\end{observation}
Note that $\opp(\square)^+=\square$ if and only if $\square\in\{(1,1),(2,k)\}$.
\medskip

In the next few definitions, we develop a larger vocabulary which will be used to discuss the movement of the frogs.

\begin{defn}[Clockwise path]\label{cpath}
    A sequence of squares $(\square_1,\dots,\square_\ell)$ coming from $[2]\times[k]$ is called a \emph{clockwise path} if $\ell\leq 2k$ and $\square_{i+1}=\square_{i}^+$ for each $i\in[\ell-1]$.
    For squares $\square_1,\square_2\in[2]\times[k]$, the (unique) clockwise path starting at $\square_1$ and ending at $\square_2$ is denoted by $I[\square_1,\square_2]$.
    We abbreviate $I(\square_1,\square_2]\eqdef I[\square_1^+,\square_2]$.
\end{defn}

Naturally, we can extend the rotation map $\rot$ to act on clockwise paths by applying it to each coordinate of the sequence.
In this way, $\rot\bigl(I[\square_1,\square_2]\bigr)=I[\rot(\square_1),\rot(\square_2)]$.

While $I[\square_1,\square_2]$ and $I(\square_1,\square_2]$ are technically sequences, we will usually treat them as sets of their elements for the sake of brevity.
Note that neither of these sets are ever empty.
It is also important to notice that, as sets, $I(\square_1,\square_2]\neq I[\square_1,\square_2]-\square_1$ if $\square_1=\square_2$.
\medskip

Now, imagine walking along a clockwise path, each of whose squares contain a frog.
When you encounter a new frog, you must place a hat upon its head; however, you must additionally ensure that the result is still a hatted-frog arrangement.
Therefore, when you visit a new frog and the opposite square is occupied by a hatted frog, you must steal the hat from that other frog so that your new froggie acquaintance can don a hat.
This idea is formalized in the following definition.

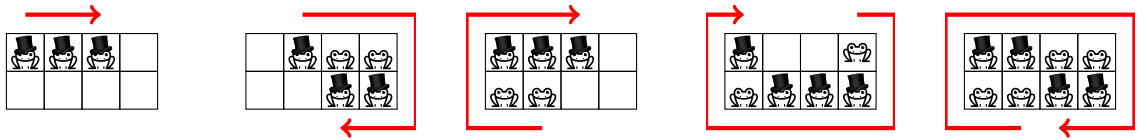
\begin{figure}[ht]
    \begin{center}
        \begin{tikzpicture}
            \clip (-0.5,-1) rectangle (2,1.5);
            \foreach \x in {0,0.5,1,1.5} {
                \draw ({\x-0.25},-0.25) rectangle (\x+0.25,0.25);
                \draw ({\x-0.25},0.25) rectangle (\x+0.25,0.75);
            }

            \node (f1) at (0,0.5) {$\froggiehat$};
            \node (f2) at (0.5,0.5) {$\froggiehat$};
            \node (f3) at (1,0.5) {$\froggiehat$};

            \draw[red, ultra thick, ->] (0,1)--(1,1);
        \end{tikzpicture}\hfil
        \begin{tikzpicture}
            \clip (-0.5,-1) rectangle (2,1.5);
            \foreach \x in {0,0.5,1,1.5} {
                \draw ({\x-0.25},-0.25) rectangle (\x+0.25,0.25);
                \draw ({\x-0.25},0.25) rectangle (\x+0.25,0.75);
            }

            \node (f2) at (0.5,0.5) {$\froggiehat$};
            \node (f3) at (1,0.5) {$\froggieb$};
            \node (f4) at (1.5,0.5) {$\froggieb$};
            \node (f5) at (1.5,0) {$\froggiehat$};
            \node (f6) at (1,0) {$\froggiehat$};

            \draw[red, ultra thick, ->] (0.5,1) -- (2,1) -- (2,-0.5) -- (1,-0.5);
        \end{tikzpicture}\hfil
        \begin{tikzpicture}
            \clip (-0.5,-1) rectangle (2,1.5);
            \foreach \x in {0,0.5,1,1.5} {
                \draw ({\x-0.25},-0.25) rectangle (\x+0.25,0.25);
                \draw ({\x-0.25},0.25) rectangle (\x+0.25,0.75);
            }

            \node (f1) at (0,0.5) {$\froggiehat$};
            \node (f2) at (0.5,0.5) {$\froggiehat$};
            \node (f3) at (1,0.5) {$\froggiehat$};
            \node (f7) at (0.5,0) {$\froggieb$};
            \node (f8) at (0,0) {$\froggieb$};

            \draw[red, ultra thick, ->] (0.5,-0.5) -- (-0.5,-0.5) -- (-0.5,1) -- (1,1);
        \end{tikzpicture}\hfil
        \begin{tikzpicture}
            \clip (-0.5,-1) rectangle (2,1.5);
            \foreach \x in {0,0.5,1,1.5} {
                \draw ({\x-0.25},-0.25) rectangle (\x+0.25,0.25);
                \draw ({\x-0.25},0.25) rectangle (\x+0.25,0.75);
            }

            \node (f1) at (0,0.5) {$\froggiehat$};
            \node (f2) at (0,0) {$\froggieb$};
            \node (f3) at (0.5,0) {$\froggiehat$};
            \node (f4) at (1,0) {$\froggiehat$};
            \node (f5) at (1.5,0) {$\froggiehat$};
            \node (f6) at (1.5,0.5) {$\froggie$};

            \draw[red, ultra thick, ->] (1.5,1) -- (2,1) -- (2,-0.5) -- (-0.5,-0.5) -- (-0.5,1) -- (0,1);
        \end{tikzpicture}\hfil
        \begin{tikzpicture}
            \clip (-0.5,-1) rectangle (2,1.5);
            \foreach \x in {0,0.5,1,1.5} {
                \draw ({\x-0.25},-0.25) rectangle (\x+0.25,0.25);
                \draw ({\x-0.25},0.25) rectangle (\x+0.25,0.75);
            }

            \node (f1) at (0,0.5) {$\froggiehat$};
            \node (f2) at (0,0) {$\froggieb$};
            \node (f3) at (0.5,0) {$\froggieb$};
            \node (f4) at (1,0) {$\froggiehat$};
            \node (f5) at (1.5,0) {$\froggiehat$};
            \node (f6) at (1.5,0.5) {$\froggieb$};
            \node (f7) at (0.5,0.5) {$\froggiehat$};
            \node (f8) at (1,0.5) {$\froggieb$};

            \draw[red, ultra thick, ->] (0.5,-0.5) -- (-0.5,-0.5) -- (-0.5,1) -- (2,1) -- (2,-0.5) -- (1,-0.5);
        \end{tikzpicture}
    \end{center}
    \caption{\label{fig:alignment}
        Examples of \hyperref[alignment]{alignment}.
        The arrow indicates a \hyperref[cpath]{clockwise path} aligning with the arrangement, though this is not the unique clockwise path with this property.
        In each example, we draw only the frogs within the clockwise path, though there could be other frogs in the arrangement.
    }
\end{figure}

\begin{defn}[Alignment]\label{alignment}
    A hatted-frog arrangement $(F,H)\in\hattedFrogs km$ is said to \emph{align} with a clockwise path $(\square_1,\dots,\square_\ell)$ in $[2]\times[k]$ if
    \begin{itemize}
        \item $\square_i\in F$ for all $i\in[\ell]$, and
        \item $\square_i\in H$ if and only if $\column(\square_i)\cap\{\square_{i+1},\dots,\square_\ell\}=\varnothing$.
    \end{itemize}
\end{defn}
See \Cref{fig:alignment} for examples of alignment.
Observe that if $\hat F$ aligns with $(\square_1,\dots,\square_\ell)$, then $\hat F$ aligns also with $(\square_i,\dots,\square_\ell)$ for every $i\in[\ell]$.

\begin{defn}[Missing squares]
    For a subset $F\subseteq[2]\times[k]$ and a square $\square\in[2]\times[k]$, write $\eb{F}(\square)$ to denote the first square preceding $\square$ in the clockwise order which is \emph{not} a member of $F$.
    If either $\square\notin F$ or $F=[2]\times[k]$, then we define $\eb{F}(\square)=\square$.
\end{defn}
Observe that if $\eb{F}(\square)\neq\square$, then $\eb{F}(\square)\notin F$ while $I(\eb{F}(\square),\square]\subseteq F$.
Additionally, notice that $\rot\bigl(\eb{F}(\square)\bigr)=\eb{\rot(F)}\bigl(\rot(\square)\bigr)$.
\medskip


The following lemma essentially states that alignment is maintained throughout the intermediate steps in the hatted-frog process.
This can be seen as an extension of \Cref{nextStep}.
\begin{lemma}\label{alignsNext}
    Fix a hatted-frog arrangement $(F,H)\in\hattedFrogs km$ and a square $\square\in[2]\times[k]$.
    Suppose that either
    \begin{itemize}
        \item $F\cap\column(\square^-)=\varnothing$, or
        \item $(F,H)$ aligns with $I(\eb{F}(\square^-),\square^-]$.
    \end{itemize}
    Then $\bigl(F+\square,\  H-\opp(\square)+\square\bigr)$ is a hatted-frog arrangement which aligns with $I(\eb{F+\square}(\square),\square]$.
\end{lemma}
\begin{proof}
    Set $(F',H')=\bigl(F+\square,\  H-\opp(\square)+\square\bigr)$.

    Note that if $(F,H)$ aligns with $I(\eb{F}(\square^-),\square^-]$, then $\square^-\in H$.
    Therefore, either $\square^-\in H$ or $F\cap\column(\square^-)=\varnothing$, and so \Cref{nextStep} tells us that $(F',H')$ is also a hatted-frog arrangement.

    Next, if $F\cap\column(\square^-)=\varnothing$, then $\eb{F'}(\square)=\square^-$ and so $I(\eb{F'}(\square),\square]=(\square)$ clearly aligns with $(F',H')$ since $\square\in H'$.

    Finally, suppose that $(F,H)$ aligns with $I(\eb{F}(\square^-),\square^-]$.
    Write $I(\eb{F}(\square^-),\square^-]=(\square_1,\dots,\square_\ell)$, so $\square_1=\eb{F}(\square^-)^+$ and $\square_\ell=\square^-$.

    Suppose first that $F+\square\neq[2]\times[k]$ or $\square\notin F$.
    It is not too difficult to see that $\eb{F+\square}(\square)=\eb{F}(\square^-)$ in this situation.
    Therefore, we have $I(\eb{F+\square}(\square),\square]=(\square_1,\dots,\square_\ell,\square)$, which aligns with $(F',H')$ since $\square\in H'$ and $(F',H')$ and $(F,H)$ agree everywhere beyond $\column(\square)$.

    Otherwise, we must have $F=[2]\times[k]$ and so $F+\square=F$.
    Thus, $\eb{F}(\square^-)=\square^-$ and $\eb{F+\square}(\square)=\square$.
    Therefore, $I(\eb{F+\square}(\square),\square]=(\square_2,\dots,\square_\ell,\square_1=\square)$.
    This path aligns with $(F',H')$ for the same reason as in the previous paragraph.
\end{proof}

\subsubsection{Crowned-frog arrangements}
Recalling the proof of \Cref{tohats}, we can identify an intermediate step in the hatted-frog process with a triple $(F,H,A)$ where $F\subseteq[2]\times[k]$ represents the positions of the un-agitated frogs, $H\subseteq F$ represents the positions of the un-agitated hatted frogs and $A\subseteq[2]\times[k]$ represents the positions of the agitated frogs.
Unfortunately, tracking the evolution of the triple $(F,H,A)$ will not be quite enough for our purposes.
As discussed at the beginning of \Cref{sec:stationary}, we also keep track of a special hatted-frog, which we call the \emph{crowned frog} and denote by $\froggiecrown$; after all, a crown is quite a special hat!

\begin{figure}
    \begin{subfigure}{\textwidth}
        \begin{center}
            \begin{tikzpicture}[baseline = (center)]
                \clip (-0.5,-1) rectangle (2,1.5);
                \foreach \x in {0,0.5,1,1.5} {
                    \draw ({\x-0.25},-0.25) rectangle (\x+0.25,0.25);
                    \draw ({\x-0.25},0.25) rectangle (\x+0.25,0.75);
                }

                \node (f1) at (0,0.5) {$\froggiehat$};
                \node (f2) at (0.5,0.5) {$\froggiehat$};
                \node (f3) at (1.5,0.5) {$\froggieb$};
                \node (f4) at (0,0) {$\froggieb$};
                \node (f3) at (1,0) {$\froggiehat$};
                \node (f4) at (1.5,0) {$\froggiehat$};
            \end{tikzpicture}
            $\xmapsto{\hyperref[poke_defn]{\poke{2}}}$
            \begin{tikzpicture}[baseline = (center)]
                \clip (-0.5,-1) rectangle (2,1.5);
                \foreach \x in {0,0.5,1,1.5} {
                    \draw ({\x-0.25},-0.25) rectangle (\x+0.25,0.25);
                    \draw ({\x-0.25},0.25) rectangle (\x+0.25,0.75);
                }

                \node (f1) at (0,0.5) {$\froggiehat$};
                \node (f2) at (0.5,1) {$\froggiecrown$};
                \node (f3) at (1.5,0.5) {$\froggieb$};
                \node (f4) at (0,0) {$\froggieb$};
                \node (f3) at (1,0) {$\froggiehat$};
                \node (f4) at (1.5,0) {$\froggiehat$};
            \end{tikzpicture}
            $\refmapsto{rule:cEH}$
            \begin{tikzpicture}[baseline = (center)]
                \clip (-0.5,-1) rectangle (2,1.5);
                \foreach \x in {0,0.5,1,1.5} {
                    \draw ({\x-0.25},-0.25) rectangle (\x+0.25,0.25);
                    \draw ({\x-0.25},0.25) rectangle (\x+0.25,0.75);
                }

                \node (f1) at (0,0.5) {$\froggiehat$};
                \node (f2) at (1,0.5) {$\froggiecrown$};
                \node (f3) at (1.5,0.5) {$\froggieb$};
                \node (f4) at (0,0) {$\froggieb$};
                \node (f3) at (1,0) {$\froggieb$};
                \node (f4) at (1.5,0) {$\froggiehat$};
            \end{tikzpicture}
        \end{center}
        \sfcaption
        \caption{Compare to \Cref{fig:hattedProcess:a}.}
    \end{subfigure}
    \sfbetween

    \begin{subfigure}{\textwidth}
        \begin{center}
            \begin{tikzpicture}[baseline = (center)]
                \clip (-0.5,-1) rectangle (2,1.5);
                \foreach \x in {0,0.5,1,1.5} {
                    \draw ({\x-0.25},-0.25) rectangle (\x+0.25,0.25);
                    \draw ({\x-0.25},0.25) rectangle (\x+0.25,0.75);
                }

                \node (f1) at (0.5,0.5) {$\froggiehat$};
                \node (f2) at (1.5,0.5) {$\froggieb$};
                \node (f3) at (1,0) {$\froggiehat$};
                \node (f4) at (1.5,0) {$\froggiehat$};
            \end{tikzpicture}
            $\xmapsto{\hyperref[poke_defn]{\poke{4}}}$
            \begin{tikzpicture}[baseline = (center)]
                \clip (-0.5,-1) rectangle (2,1.5);
                \foreach \x in {0,0.5,1,1.5} {
                    \draw ({\x-0.25},-0.25) rectangle (\x+0.25,0.25);
                    \draw ({\x-0.25},0.25) rectangle (\x+0.25,0.75);
                }

                \node (f1) at (0.5,0.5) {$\froggiehat$};
                \node (f2) at (1.5,1) {$\froggie$};
                \node (f3) at (1,0) {$\froggiehat$};
                \node (f4) at (1.5,-0.5) {$\froggiecrown$};
            \end{tikzpicture}
            $\refmapsto{rule:fE}$
            \begin{tikzpicture}[baseline = (center)]
                \clip (-0.5,-1) rectangle (2,1.5);
                \foreach \x in {0,0.5,1,1.5} {
                    \draw ({\x-0.25},-0.25) rectangle (\x+0.25,0.25);
                    \draw ({\x-0.25},0.25) rectangle (\x+0.25,0.75);
                }

                \node (f1) at (0.5,0.5) {$\froggiehat$};
                \node (f2) at (1.5,0) {$\froggiehat$};
                \node (f3) at (1,0) {$\froggiehat$};
                \node (f4) at (1.5,-0.5) {$\froggiecrown$};
            \end{tikzpicture}
            $\refmapsto{rule:cH}$
            \begin{tikzpicture}[baseline = (center)]
                \clip (-0.5,-1) rectangle (2,1.5);
                \foreach \x in {0,0.5,1,1.5} {
                    \draw ({\x-0.25},-0.25) rectangle (\x+0.25,0.25);
                    \draw ({\x-0.25},0.25) rectangle (\x+0.25,0.75);
                }

                \node (f1) at (0.5,0.5) {$\froggiehat$};
                \node (f2) at (1.5,0) {$\froggiehat$};
                \node (f3) at (1,-0.5) {$\froggiecrown$};
                \node (f4) at (1,0) {$\froggiehat$};
            \end{tikzpicture}
            $\refmapsto{rule:cEH}$
            \begin{tikzpicture}[baseline = (center)]
                \clip (-0.5,-1) rectangle (2,1.5);
                \foreach \x in {0,0.5,1,1.5} {
                    \draw ({\x-0.25},-0.25) rectangle (\x+0.25,0.25);
                    \draw ({\x-0.25},0.25) rectangle (\x+0.25,0.75);
                }

                \node (f1) at (0.5,0.5) {$\froggieb$};
                \node (f2) at (1.5,0) {$\froggiehat$};
                \node (f3) at (0.5,0) {$\froggiecrown$};
                \node (f4) at (1,0) {$\froggiehat$};
            \end{tikzpicture}
        \end{center}
        \sfcaption
        \caption{Compare to \Cref{fig:hattedProcess:b}.}
    \end{subfigure}
    \sfbetween

    \begin{subfigure}{\textwidth}
        \begin{center}
            \begin{tikzpicture}[baseline = (center)]
                \clip (-0.5,-1) rectangle (2,1.5);
                \foreach \x in {0,0.5,1,1.5} {
                    \draw ({\x-0.25},-0.25) rectangle (\x+0.25,0.25);
                    \draw ({\x-0.25},0.25) rectangle (\x+0.25,0.75);
                }

                \node (f1) at (0,0.5) {$\froggiehat$};
                \node (f2) at (0.5,0.5) {$\froggieb$};
                \node (f3) at (1,0.5) {$\froggieb$};
                \node (f4) at (0,0) {$\froggieb$};
                \node (f5) at (0.5,0) {$\froggiehat$};
                \node (f6) at (1,0) {$\froggiehat$};
                \node (f7) at (1.5,0) {$\froggiehat$};
            \end{tikzpicture}
            $\xmapsto{\hyperref[poke_defn]{\poke{2}}}$
            \begin{tikzpicture}[baseline = (center)]
                \clip (-0.5,-1) rectangle (2,1.5);
                \foreach \x in {0,0.5,1,1.5} {
                    \draw ({\x-0.25},-0.25) rectangle (\x+0.25,0.25);
                    \draw ({\x-0.25},0.25) rectangle (\x+0.25,0.75);
                }

                \node (f1) at (0,0.5) {$\froggiehat$};
                \node (f2) at (0.5,1) {$\froggie$};
                \node (f3) at (1,0.5) {$\froggieb$};
                \node (f4) at (0,0) {$\froggieb$};
                \node (f5) at (0.5,-0.5) {$\froggiecrown$};
                \node (f6) at (1,0) {$\froggiehat$};
                \node (f7) at (1.5,0) {$\froggiehat$};
            \end{tikzpicture}
            $\refmapsto{rule:fFH}$
            \begin{tikzpicture}[baseline = (center)]
                \clip (-0.5,-1) rectangle (2,1.5);
                \foreach \x in {0,0.5,1,1.5} {
                    \draw ({\x-0.25},-0.25) rectangle (\x+0.25,0.25);
                    \draw ({\x-0.25},0.25) rectangle (\x+0.25,0.75);
                }

                \node (f1) at (0,0.5) {$\froggiehat$};
                \node (f2) at (1,0.5) {$\froggiehat$};
                \node (f3) at (1,1) {$\froggie$};
                \node (f4) at (0,0) {$\froggieb$};
                \node (f5) at (0.5,-0.5) {$\froggiecrown$};
                \node (f6) at (1,0) {$\froggieb$};
                \node (f7) at (1.5,0) {$\froggiehat$};
            \end{tikzpicture}
            $\refmapsto{rule:fEH}$
            \begin{tikzpicture}[baseline = (center)]
                \clip (-0.5,-1) rectangle (2,1.5);
                \foreach \x in {0,0.5,1,1.5} {
                    \draw ({\x-0.25},-0.25) rectangle (\x+0.25,0.25);
                    \draw ({\x-0.25},0.25) rectangle (\x+0.25,0.75);
                }

                \node (f1) at (0,0.5) {$\froggiehat$};
                \node (f2) at (1,0.5) {$\froggiehat$};
                \node (f3) at (1.5,0.5) {$\froggiehat$};
                \node (f4) at (0,0) {$\froggieb$};
                \node (f5) at (0.5,-0.5) {$\froggiecrown$};
                \node (f6) at (1,0) {$\froggieb$};
                \node (f7) at (1.5,0) {$\froggieb$};
            \end{tikzpicture}
            $\refmapsto{rule:cFH}$
            \begin{tikzpicture}[baseline = (center)]
                \clip (-0.5,-1) rectangle (2,1.5);
                \foreach \x in {0,0.5,1,1.5} {
                    \draw ({\x-0.25},-0.25) rectangle (\x+0.25,0.25);
                    \draw ({\x-0.25},0.25) rectangle (\x+0.25,0.75);
                }

                \node (f1) at (0,0.5) {$\froggieb$};
                \node (f2) at (1,0.5) {$\froggiehat$};
                \node (f3) at (1.5,0.5) {$\froggiehat$};
                \node (f4) at (0,-0.5) {$\froggie$};
                \node (f5) at (0,0) {$\froggiecrown$};
                \node (f6) at (1,0) {$\froggieb$};
                \node (f7) at (1.5,0) {$\froggieb$};
            \end{tikzpicture}
            $\refmapsto{rule:fFC}$
            \begin{tikzpicture}[baseline = (center)]
                \clip (-0.5,-1) rectangle (2,1.5);
                \foreach \x in {0,0.5,1,1.5} {
                    \draw ({\x-0.25},-0.25) rectangle (\x+0.25,0.25);
                    \draw ({\x-0.25},0.25) rectangle (\x+0.25,0.75);
                }

                \node (f1) at (0,1) {$\froggiecrown$};
                \node (f2) at (1,0.5) {$\froggiehat$};
                \node (f3) at (1.5,0.5) {$\froggiehat$};
                \node (f4) at (0,0.5) {$\froggiehat$};
                \node (f5) at (0,0) {$\froggieb$};
                \node (f6) at (1,0) {$\froggieb$};
                \node (f7) at (1.5,0) {$\froggieb$};
            \end{tikzpicture}
            $\refmapsto{rule:cE}$
            \begin{tikzpicture}[baseline = (center)]
                \clip (-0.5,-1) rectangle (2,1.5);
                \foreach \x in {0,0.5,1,1.5} {
                    \draw ({\x-0.25},-0.25) rectangle (\x+0.25,0.25);
                    \draw ({\x-0.25},0.25) rectangle (\x+0.25,0.75);
                }

                \node (f1) at (0.5,0.5) {$\froggiecrown$};
                \node (f2) at (1,0.5) {$\froggiehat$};
                \node (f3) at (1.5,0.5) {$\froggiehat$};
                \node (f4) at (0,0.5) {$\froggiehat$};
                \node (f5) at (0,0) {$\froggieb$};
                \node (f6) at (1,0) {$\froggieb$};
                \node (f7) at (1.5,0) {$\froggieb$};
            \end{tikzpicture}
        \end{center}
        \sfcaption
        \caption{Compare to \Cref{fig:hattedProcess:c}.}
    \end{subfigure}
    \sfbetween

    \begin{subfigure}{\textwidth}
        \begin{center}
            \begin{tikzpicture}[baseline = (center)]
                \clip (-0.5,-1) rectangle (2,1.5);
                \foreach \x in {0,0.5,1,1.5} {
                    \draw ({\x-0.25},-0.25) rectangle (\x+0.25,0.25);
                    \draw ({\x-0.25},0.25) rectangle (\x+0.25,0.75);
                }

                \node (f1) at (0,0.5) {$\froggiehat$};
                \node (f2) at (0.5,0.5) {$\froggiehat$};
                \node (f3) at (1,0.5) {$\froggieb$};
                \node (f4) at (0,0) {$\froggieb$};
                \node (f5) at (0.5,0) {$\froggieb$};
                \node (f6) at (1,0) {$\froggiehat$};
                \node (f7) at (1.5,0) {$\froggiehat$};
            \end{tikzpicture}
            $\xmapsto{\hyperref[poke_defn]{\poke{3}}}$
            \begin{tikzpicture}[baseline = (center)]
                \clip (-0.5,-1) rectangle (2,1.5);
                \foreach \x in {0,0.5,1,1.5} {
                    \draw ({\x-0.25},-0.25) rectangle (\x+0.25,0.25);
                    \draw ({\x-0.25},0.25) rectangle (\x+0.25,0.75);
                }

                \node (f1) at (0,0.5) {$\froggiehat$};
                \node (f2) at (0.5,0.5) {$\froggiehat$};
                \node (f3) at (1,1) {$\froggie$};
                \node (f4) at (0,0) {$\froggieb$};
                \node (f5) at (0.5,0) {$\froggieb$};
                \node (f6) at (1,-0.5) {$\froggiecrown$};
                \node (f7) at (1.5,0) {$\froggiehat$};
            \end{tikzpicture}
            $\refmapsto{rule:fEH}$
            \begin{tikzpicture}[baseline = (center)]
                \clip (-0.5,-1) rectangle (2,1.5);
                \foreach \x in {0,0.5,1,1.5} {
                    \draw ({\x-0.25},-0.25) rectangle (\x+0.25,0.25);
                    \draw ({\x-0.25},0.25) rectangle (\x+0.25,0.75);
                }

                \node (f1) at (0,0.5) {$\froggiehat$};
                \node (f2) at (0.5,0.5) {$\froggiehat$};
                \node (f3) at (1.5,0.5) {$\froggiehat$};
                \node (f4) at (0,0) {$\froggieb$};
                \node (f5) at (0.5,0) {$\froggieb$};
                \node (f6) at (1,-0.5) {$\froggiecrown$};
                \node (f7) at (1.5,0) {$\froggieb$};
            \end{tikzpicture}
            $\refmapsto{rule:cFH}$
            \begin{tikzpicture}[baseline = (center)]
                \clip (-0.5,-1) rectangle (2,1.5);
                \foreach \x in {0,0.5,1,1.5} {
                    \draw ({\x-0.25},-0.25) rectangle (\x+0.25,0.25);
                    \draw ({\x-0.25},0.25) rectangle (\x+0.25,0.75);
                }

                \node (f1) at (0,0.5) {$\froggiehat$};
                \node (f2) at (0.5,0.5) {$\froggieb$};
                \node (f3) at (1.5,0.5) {$\froggiehat$};
                \node (f4) at (0,0) {$\froggieb$};
                \node (f5) at (0.5,-0.5) {$\froggie$};
                \node (f6) at (0.5,0) {$\froggiecrown$};
                \node (f7) at (1.5,0) {$\froggieb$};
            \end{tikzpicture}
            $\refmapsto{rule:fFH}$
            \begin{tikzpicture}[baseline = (center)]
                \clip (-0.5,-1) rectangle (2,1.5);
                \foreach \x in {0,0.5,1,1.5} {
                    \draw ({\x-0.25},-0.25) rectangle (\x+0.25,0.25);
                    \draw ({\x-0.25},0.25) rectangle (\x+0.25,0.75);
                }

                \node (f1) at (0,0.5) {$\froggieb$};
                \node (f2) at (0.5,0.5) {$\froggieb$};
                \node (f3) at (1.5,0.5) {$\froggiehat$};
                \node (f4) at (0,-0.5) {$\froggie$};
                \node (f5) at (0,0) {$\froggiehat$};
                \node (f6) at (0.5,0) {$\froggiecrown$};
                \node (f7) at (1.5,0) {$\froggieb$};
            \end{tikzpicture}
            $\refmapsto{rule:fFH}$
            \begin{tikzpicture}[baseline = (center)]
                \clip (-0.5,-1) rectangle (2,1.5);
                \foreach \x in {0,0.5,1,1.5} {
                    \draw ({\x-0.25},-0.25) rectangle (\x+0.25,0.25);
                    \draw ({\x-0.25},0.25) rectangle (\x+0.25,0.75);
                }

                \node (f1) at (0,1) {$\froggie$};
                \node (f2) at (0.5,0.5) {$\froggieb$};
                \node (f3) at (1.5,0.5) {$\froggiehat$};
                \node (f4) at (0,0.5) {$\froggiehat$};
                \node (f5) at (0,0) {$\froggieb$};
                \node (f6) at (0.5,0) {$\froggiecrown$};
                \node (f7) at (1.5,0) {$\froggieb$};
            \end{tikzpicture}
            $\refmapsto{rule:fFC}$
            \begin{tikzpicture}[baseline = (center)]
                \clip (-0.5,-1) rectangle (2,1.5);
                \foreach \x in {0,0.5,1,1.5} {
                    \draw ({\x-0.25},-0.25) rectangle (\x+0.25,0.25);
                    \draw ({\x-0.25},0.25) rectangle (\x+0.25,0.75);
                }

                \node (f1) at (0.5,0.5) {$\froggiehat$};
                \node (f2) at (0.5,1) {$\froggiecrown$};
                \node (f3) at (1.5,0.5) {$\froggiehat$};
                \node (f4) at (0,0.5) {$\froggiehat$};
                \node (f5) at (0,0) {$\froggieb$};
                \node (f6) at (0.5,0) {$\froggieb$};
                \node (f7) at (1.5,0) {$\froggieb$};
            \end{tikzpicture}
            $\refmapsto{rule:cE}$
            \begin{tikzpicture}[baseline = (center)]
                \clip (-0.5,-1) rectangle (2,1.5);
                \foreach \x in {0,0.5,1,1.5} {
                    \draw ({\x-0.25},-0.25) rectangle (\x+0.25,0.25);
                    \draw ({\x-0.25},0.25) rectangle (\x+0.25,0.75);
                }

                \node (f1) at (0.5,0.5) {$\froggiehat$};
                \node (f2) at (1,0.5) {$\froggiecrown$};
                \node (f3) at (1.5,0.5) {$\froggiehat$};
                \node (f4) at (0,0.5) {$\froggiehat$};
                \node (f5) at (0,0) {$\froggieb$};
                \node (f6) at (0.5,0) {$\froggieb$};
                \node (f7) at (1.5,0) {$\froggieb$};
            \end{tikzpicture}
        \end{center}
        \sfcaption
        \caption{Compare to \Cref{fig:hattedProcess:d}.}
    \end{subfigure}
    \caption{\label{fig:crownedProcess}
        Examples of \Cref{crown-coupling}.
        An agitated frog is drawn just above/below the square that it occupies.
        The examples shown here correspond to the examples of the hatted-frog process shown in \Cref{fig:hattedProcess}.
    }
\end{figure}
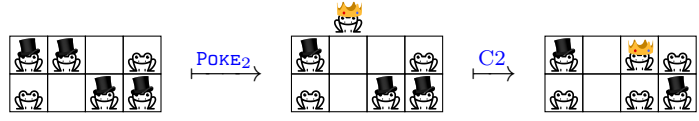
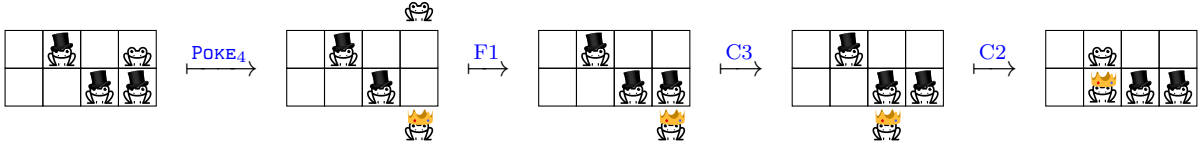
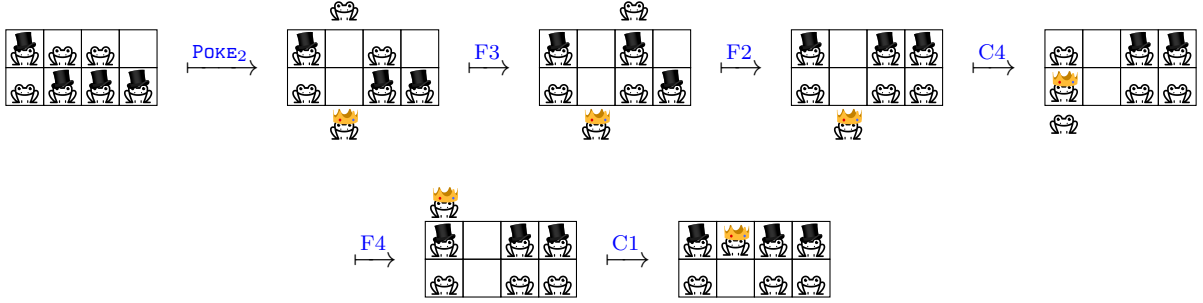
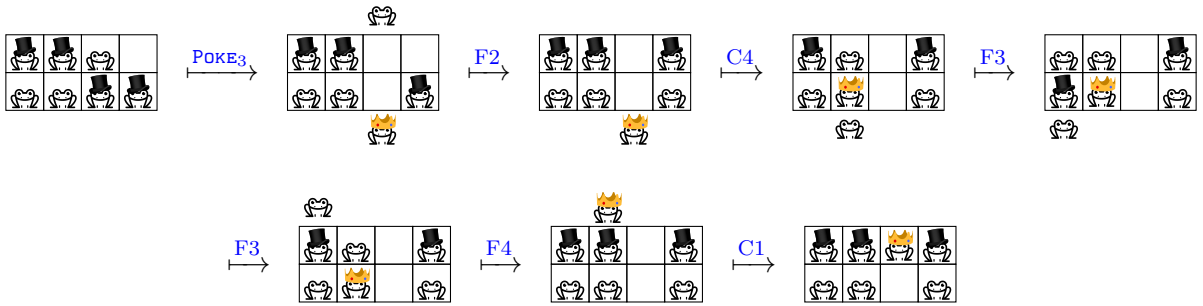

The following definition is technical; the reader may find it helpful to peek ahead to \Cref{fig:crownedProcess} to see examples, even though much of that figure will not yet make sense.

\begin{defn}[Crowned-frog arrangement]
    Fix positive integers $k,m$.
    The set of \emph{$m$-crowned-frog arrangements} on $[2]\times[k]$ consists of all quintuples of the form $(F,H,A,\froggiecrown,x)$ where $F,H,A\subseteq[2]\times[k]$, and $\froggiecrown\in [2]\times[k]$, and $x\in\{\agitated,\settled\}$, satisfying the following properties:
    \begin{properties}
        \item\label{cf:hatted} $(F,H)$ is a hatted-frog arrangement.
        \item\label{cf:sum} $\abs A\leq 2$ and $\abs F+\abs A=m$.
        \item\label{cf:whereCrown} $\froggiecrown\in H\cup A$, where
            \[
                x=\agitated\implies\froggiecrown\in A,\qquad\text{and}\qquad x=\settled\implies\froggiecrown\in H.
            \]
        \item\label{cf:hatOrEmpty} For each $\froggie\in A$, either $F\cap\column(\froggie)=\varnothing$ or $(F,H)$ aligns with $I(\eb{F}(\froggie),\froggie]$.
        \item\label{cf:one} If $A=\{\froggie\}$ and $x=\settled$, then $\froggie^+\notin H$ and $\column(\square)\subseteq F$ for every $\square\in I[\froggiecrown,\froggie]$.
        \item\label{cf:two} If $\abs A=2$, then
            \begin{properties}
                \item\label{cf:two:empty} $F\cap\column(\froggiecrown)=\varnothing$ (which implies that $\froggiecrown\in A\setminus H$ and $x=\agitated$), and
                \item\label{cf:two:other} If $A=\{\froggiecrown,\froggie\}$, then $\froggie^+\notin H$ and $\column(\square)\subseteq F$ for every $\square\in I[\opp(\froggiecrown),\froggie]-\opp(\froggiecrown)$.
            \end{properties}
                    \item\label{cf:agitated} If $x=\agitated$ and $\eb{F}(\froggiecrown)=\opp(\froggiecrown^+)$, then $F\cap\column(\froggiecrown^+)=\varnothing$.
    \end{properties}
    The set of all $m$-crowned-frog arrangements on $[2]\times[k]$ is denoted by $\crownedFrogs km$.
\end{defn}

As suggested earlier, the set $F$ represents the positions of the un-agitated frogs, the set $H$ represents the positions of the un-agitated frogs wearing a hat and the set $A$ represents the positions of the agitated frogs.
Furthermore, $\froggiecrown$ marks the position of a special crowned frog (which is either a hatted frog or an agitated frog) and $x\in\{\agitated,\settled\}$ indicates whether or not the crowned frog is agitated (which is necessary only when $\froggiecrown\in A\cap H$ due to \cref{cf:whereCrown}).

We extend the rotation map $\rot$ to act on $\crownedFrogs km$ by applying it pointwise to the first four coordinates and leaving the fifth unchanged.
It is not difficult to see that $\crownedFrogs km$ is invariant under $\rot$.
\medskip

We next define two special subsets of $\crownedFrogs km$ which bridge the gap between crowned-frog arrangements and hatted-frog arrangements.
\begin{defn}[Starting arrangements]
    The set $\crownedFrogsStart km$ consists of all crowned-frog arrangements of the form $(F,H,A,\froggiecrown,\agitated)\in\crownedFrogs km$ where $A\subseteq\column(\froggiecrown)$, and $F\cap\column(\froggiecrown)=\varnothing$, and ${(F\sqcup A,H+\froggiecrown)}\in\hattedFrogs km$.
\end{defn}
\begin{defn}[Ending arrangements]
    The set $\crownedFrogsEnd km$ consists of all crowned-frog arrangements of the form $(F,H,\varnothing,\froggiecrown,\settled)\in\crownedFrogs km$.

\end{defn}
We now justify the relationship between starting/ending arrangements and hatted-frog arrangements.
\medskip

For each $c\in[k]$, define a map
\[
    \poke c\colon\hattedFrogs km\to\crownedFrogs km\sqcup\hattedFrogs km,
\]
which ``pokes the frogs in column $c$''.
For $\hat F=(F,H)\in\hattedFrogs km$, set $A=F\cap\bigl([2]\times\{c\}\bigr)$:
\begin{itemize}
    \item If $A=\varnothing$, then $\poke c(\hat F)\eqdef\hat F$.
    \item\label{poke_defn} If $A\neq\varnothing$, set $\{\froggiecrown\}=H\cap A$.
        Then $\poke c(\hat F)\eqdef(F\setminus A,\ H\setminus A,\ A,\ \froggiecrown,\ \agitated)$.
\end{itemize}
In other words, the frogs in column $c$ become agitated and the frog in column $c$ which had a hat now dons a crown.

\begin{observation}\label{injective-poke}
    For each $\hat F\in\hattedFrogs km$ and each $c\in[k]$, $\poke{c}(\hat F)\in\crownedFrogsStart km\sqcup\hattedFrogs km$.
    Furthermore, for each $\crown F\in\crownedFrogsStart km$, there is a \emph{unique} pair $(\hat F,c)\in\hattedFrogs km\times[k]$ for which $\crown F=\poke{c}(\hat F)$.
\end{observation}

Turning our attention to $\crownedFrogsEnd km$, there is a projection which simply forgets about the crowned frog:
\begin{align*}
    \dethrone &\colon \crownedFrogsEnd km\to\hattedFrogs km\\
              &\colon (F,H,\varnothing,\froggiecrown,\settled)\mapsto (F,H).
\end{align*}
\begin{observation}\label{ending-select}
    The elements of $\crownedFrogsEnd km$ are created by taking an element of $\hattedFrogs km$ and declaring one of the hatted frogs to have a crown.
    In particular, for each $\hat F=(F,H)\in\hattedFrogs km$,
    \[
        \dethrone^{-1}(\hat F)=\bigl\{(F,H,\varnothing,\froggiehat,\settled):\froggiehat\in H\bigr\}\quad\implies\quad\abs{\dethrone^{-1}(\hat F)}=\abs H.
    \]
\end{observation}

With these definitions in place, we encourage the reader to look back at the discussion toward the beginning of \Cref{sec:stationary}.
Essentially, proving \Cref{regular} amounts to finding a bijection $\Phi\colon\crownedFrogsStart km\to\crownedFrogsEnd km$ which tracks the hatted-frog process in the sense that
\[
    \hat Fc=(\dethrone\circ\,\Phi\circ\poke{c})(\hat F),
\]
whenever $\hat F$ is a hatted-frog arrangement with a frog in column $c$.

To find this bijection, we define a transition map on $\crownedFrogs km$ which represents the hopping of an agitated frog.
In the following definition, we include also pictorial representations of each of the rules to aid in understanding.
We additionally encourage the reader to look at \Cref{fig:crownedProcess} while reading the definition in order to see the transition map in action.

\begin{defn}[The transition map $\move$]\label{movedefn}
    Fix an element $\crown F=(F,H,A,\froggiecrown,x)\in\crownedFrogs km\setminus\crownedFrogsEnd km$.
    In order to define $\move(\crown F)$, we break into cases:
    \begin{rules}
        \setcounter{rulesi}{2}
        \item\label{rule:crownjump} $A=\{\froggiecrown\}$ and $x=\agitated$.
            The included pictures show only $\column(\froggiecrown^+)$ and assume that $\froggiecrown^+$ is a square in the top row of the grid.
            \begin{rules}
                \begin{minipage}{0.68\textwidth}
                \item\label{rule:cE}
                    If $F\cap\column(\froggiecrown^+)=\varnothing$, then
                    \[
                        \move(\crown F) \eqdef  \bigl(  F+\froggiecrown^+,  H+\froggiecrown^+,  \varnothing,  \froggiecrown^+,  \settled  \bigr)
                    \]
                \end{minipage}
                \begin{minipage}{0.2\textwidth}
                    \[
                        \dominoTransition{}{}{\froggiecrown}\ \mapsto\dominoTransitioned{\froggiecrown}{}{}
                    \]
                \end{minipage}
                \medskip

                \begin{minipage}{0.68\textwidth}
                \item\label{rule:cEH} If $\froggiecrown^+\notin F$ and $F\cap\column(\froggiecrown^+)\neq\varnothing$, then
                    \[
                        \move(\crown F) \eqdef \bigl( F+\froggiecrown^+,  H-\opp(\froggiecrown^+)+\froggiecrown^+, \varnothing, \froggiecrown^+, \settled \bigr)
                    \]
                \end{minipage}
                \begin{minipage}{0.2\textwidth}
                    \[
                        \dominoTransition{}{\froggiehat}{\froggiecrown}\ \mapsto \dominoTransitioned{\froggiecrown}{\froggieb}{}
                    \]
                \end{minipage}
                \medskip

                \begin{minipage}{0.68\textwidth}
                \item\label{rule:cH} If $\froggiecrown^+\in H$, then
                    \[
                        \move(\crown F) \eqdef  \bigl(  F,  H,  \{\froggiecrown^+\},  \froggiecrown^+,  \agitated  \bigr)
                    \]
                \end{minipage}
                \begin{minipage}{0.2\textwidth}
                    \[
                        \dominoTransition{\froggiehat}{?}{\froggiecrown}\ \mapsto \dominoTransitioned{\froggiehat}{?}{\froggiecrown}
                    \]
                \end{minipage}
                \medskip

                \begin{minipage}{0.68\textwidth}
                \item\label{rule:cFH} If $\froggiecrown^+\in F\setminus H$, then
                    \[
                        \move(\crown F)\eqdef\bigl(F, H-\opp(\froggiecrown^+)+\froggiecrown^+,\{\froggiecrown^+\},\froggiecrown^+, \settled\bigr)
                    \]
                \end{minipage}
                \begin{minipage}{0.2\textwidth}
                    \[
                        \dominoTransition{\froggieb}{\froggiehat}{\froggiecrown}\ \mapsto \dominoTransitioned{\froggiecrown}{\froggieb}{\froggie}
                    \]
                \end{minipage}

            \end{rules}
        \setcounter{rulesi}{5}
        \item\label{rule:uncrownjump} $\abs A=2$ or both $\abs A=1$ and $x=\settled$.
            Thus, either $A=\{\froggie\}$ or $A=\{\froggie,\froggiecrown\}$ for some $\froggie\in[2]\times[k]$.
            Note that $\froggie^+\notin H$ according to \cref{cf:one,cf:two}.
            The included pictures show only $\column(\froggie^+)$ and assume that $\froggie^+$ is a square in the top row of the grid.
            \begin{rules}
                \begin{minipage}{0.68\textwidth}
                \item\label{rule:fE} If $F\cap\column(\froggie^+)=\varnothing$, then
                    \[
                        \move(\crown F)\eqdef\bigl(F+\froggie^+,H+\froggie^+,A-\froggie,\froggiecrown,x\bigr).
                    \]
                \end{minipage}
                \begin{minipage}{0.2\textwidth}
                    \[
                        \dominoTransition{}{}{\froggie}\ \mapsto \dominoTransitioned{\froggiehat}{}{}
                    \]
                \end{minipage}
                \medskip

                \begin{minipage}{0.68\textwidth}
                \item\label{rule:fEH} If $\froggie^+\notin F$ and $F\cap\column(\froggie^+)\neq\varnothing$, then
                    \[
                        \move(\crown F)\eqdef\bigl(F+\froggie^+, H-\opp(\froggie^+)+\froggie^+,A-\froggie,\froggiecrown,x\bigr).
                    \]
                \end{minipage}
                \begin{minipage}{0.2\textwidth}
                    \[
                        \dominoTransition{}{\froggiehat}{\froggie}\ \mapsto \dominoTransitioned{\froggiehat}{\froggieb}{}
                    \]
                \end{minipage}
                \medskip


                \begin{minipage}{0.68\textwidth}
                \item\label{rule:fFH} If $\froggie^+\in F\setminus H$ and $\froggiecrown\neq\opp(\froggie^+)$, then
                    \[
                        \move(\crown F)\eqdef\bigl(F, H-\opp(\froggie^+)+\froggie^+,A-\froggie+\froggie^+,\froggiecrown,x\bigr).
                    \]
                \end{minipage}
                \begin{minipage}{0.2\textwidth}
                    \[
                        \dominoTransition{\froggieb}{\froggiehat}{\froggie}\ \mapsto \dominoTransitioned{\froggiehat}{\froggieb}{\froggie}
                    \]
                \end{minipage}
                \medskip

                \begin{minipage}{0.68\textwidth}
                \item\label{rule:fFC} If $\froggie^+\in F\setminus H$ and $\froggiecrown=\opp(\froggie^+)$, then
                    \[
                        \move(\crown F)\eqdef\bigl(F, H-\opp(\froggie^+)+\froggie^+,\{\froggie^+\},\froggie^+,\agitated\bigr).
                    \]
                \end{minipage}
                \begin{minipage}{0.2\textwidth}
                    \[
                        \dominoTransition{\froggieb}{\froggiecrown}{\froggie}\ \mapsto \dominoTransitioned{\froggiehat}{\froggieb}{\froggiecrown}
                    \]
                \end{minipage}
            \end{rules}
    \end{rules}
\end{defn}
Observe that $\move\circ\rot=\rot\circ\move$.
\medskip

At this point, we encourage the reader to skip ahead to \Cref{crown-coupling} and \Cref{fig:crownedProcess} to get a better understanding of $\move$ before diving into \Cref{well-defined}.
\medskip

We record a few observations that follow directly from the definition of $\move$.
\begin{observation}\label{moveObservations}
    Fix an element $\crown F=(F,H,A,\froggiecrown,x)\in\crownedFrogs km\setminus\crownedFrogsEnd km$.
    Define $\froggie\in[2]\times[k]$ by $\{\froggie,\froggiecrown\}=A+\froggiecrown$ (note that we could have $\froggie=\froggiecrown$).
    If $\move(\crown F)=(F_*,H_*,A_*,\froggiecrown_*,x_*)$, then
    \begin{itemize}
        \item $F_*=F+\froggie^+$ and $H_*= H-\opp(\froggie^+)+\froggie^+$, and
        \item $\froggiecrown_*\in\{\froggiecrown,\froggie^+\}$, and
        \item $\displaystyle A_*=\begin{cases}
                A-\froggie, & \text{if }\froggie^+\notin F,\\
                A-\froggie+\froggie^+, & \text{if }\froggie^+\in F.
            \end{cases}$
    \end{itemize}
\end{observation}

We now prove that $\move$ sends crowned-frog arrangements to crowned-frog arrangements.

\begin{theorem}\label{well-defined}
    Fix integers $k,m$ with $0\leq m\leq 2k$.
    If $\crown F\in\crownedFrogs km\setminus\crownedFrogsEnd km$, then $\move(\crown F)\in\crownedFrogs km\setminus\crownedFrogsStart km$.
\end{theorem}
\begin{proof}
    Write $\crown F_*=(F_*,H_*,A_*,\froggiecrown_*,x_*)=\move(\crown F)$ and $\{\froggie\}=A\setminus A_*$.
    \medskip

    We begin by showing that $\crown F_*\in\crownedFrogs km$ by verifying each of the 7 properties.
    \begin{properties}
        \item Since $\crown F \in\crownedFrogs km$, we know that $(F,H)$ is a hatted-frog arrangement (\cref{cf:hatted}) and that either $F\cap\column(\froggie)=\varnothing$ or $(F,H)$ aligns with $I(\eb F(\froggie),\froggie]$ (\cref{cf:hatOrEmpty}).
            Thus, the fact that $(F_*,H_*)$ is a hatted-frog arrangement follows by applying \Cref{alignsNext} with $\square=\froggie^+$.
        \item The fact that $\abs{A_*}\leq 2$ and $\abs{A_*}+\abs{F_*}=m$ is an immediate consequence of \Cref{moveObservations}.
        \item Suppose first that $x_*=\agitated$; thus the transition from $\crown F$ to $\crown F_*$ \emph{cannot} have taken place according to \crefOr{rule:cE,rule:cEH,rule:cFH}.
            If the transition took place according to either \crefOr{rule:cH,rule:fFC}, then $\froggiecrown_*\in A_*$ by definition.

            Otherwise, the transition took place according to one of \crefOr{rule:fE,rule:fEH,rule:fFH} and $x=x_*=\agitated$; this additionally implies that $\abs A=2$.
            Of course, in each of these rules, $\froggiecrown_*=\froggiecrown$; furthermore, $A=\{\froggie,\froggiecrown\}$ due to \cref{cf:whereCrown}.
            Since $\{\froggie\}=A\setminus A_*$, we conclude that $\froggiecrown_*=\froggiecrown\in A_*$ as needed.
            \medskip

            Suppose next that $x_*=\settled$; thus the transition from $\crown F$ to $\crown F_*$ \emph{cannot} have taken place according to \crefOr{rule:cFH,rule:fFC}.
            If the transition took place according to \crefOr{rule:cE,rule:cEH,rule:cFH}, then $\froggiecrown_*\in H_*$ by definition.

            Otherwise, the transition took place according to one of \crefOr{rule:fE,rule:fEH,rule:fFH} and $x=\settled$.
            In this case, $A=\{\froggie\}$ and $\froggiecrown_*=\froggiecrown$.
            Therefore, $\froggiecrown_*=\froggiecrown\in H$ due to \cref{cf:whereCrown}.
            Now, since $H_*= H-\opp(\froggie^+)+\froggie^+$, the only way to have $\froggiecrown_*\notin H_*$ is if $\froggiecrown_*=\opp(\froggie^+)$.
            However, we know this to not be the case since otherwise the transition would have taken place according to \cref{rule:fFC}, meaning that $x_*=\agitated$.

        \item Note that $(F,H)$ is a hatted-frog arrangement (\cref{cf:hatted}) and either $F\cap\column(\froggie)=\varnothing$ or $(F,H)$ aligns with $I(\eb{F}(\froggie),\froggie]$ (\cref{cf:hatOrEmpty}).
            We thus know that $(F_*,H_*)$ aligns with $I(\eb{F_*}(\froggie^+),\froggie^+]$ by applying \Cref{alignsNext} with $\square=\froggie^+$.

            Now, fix $\froggie_*\in A_*$; we must show either that $(F_*,H_*)$ aligns with $I(\eb{F_*}(\froggie_*),\froggie_*]$ or $F_* \cap \column(\froggie_*) = \varnothing$.
            Of course, we may suppose that $\froggie_*\neq\froggie^+$ due to the previous paragraph.
            We claim that $F_*\cap\column(\froggie_*)=\varnothing$ in this case.

            According to \Cref{moveObservations}, either $A_*=A-\froggie$ or $A_*=A-\froggie+\froggie^+$, so $\froggie_*\in A\cap A_*$ since $\froggie_*\neq\froggie^+$.
            Since $\{\froggie\}=A\setminus A_*$, we find that $\abs A=2$ and $A=\{\froggie,\froggie_*\}$.
            The fact that $\crown F$ satisfies \cref{cf:two} then forces $\froggie_*=\froggiecrown$ and $F\cap\column(\froggie_*)=\varnothing$.
            Therefore, again according to \Cref{moveObservations}, the only way that we could have $F_*\cap\column(\froggie_*)\neq\varnothing$ is if $\froggie^+\in\column(\froggie_*)$.
            Since $\froggie^+\neq\froggie_*$, this would require $\froggie^+=\opp(\froggie_*)=\opp(\froggiecrown)$.
            However, then $\froggiecrown\in I[\opp(\froggiecrown),\froggie]-\opp(\froggiecrown)$ and so \cref{cf:two:other} of $\crown F$ would imply that $\column(\froggiecrown)\subseteq F$; a contradiction.

        \item Suppose that $A_*=\{\froggie_*\}$ and $x_*=\settled$.
            We claim that $\froggie_*=\froggie^+$ in this case.
            Indeed, if this were not the case, then $A_*=A-\froggie$ (\cref{moveObservations}) and so $\abs A=2$.
            As such, $x=\agitated$ due to \cref{cf:two:empty} of $\crown F$, and so, by the definiton of $\move$, we must also have $x_*=\agitated$; a contradiction.

            Hence $\froggie_*=\froggie^+$.
            From this, we find that the transition from $\crown F$ to $\crown F_*$ must have taken place according to \crefOr{rule:cFH,rule:fFH}.
            In either case, we know that $F=F_*$, $\column(\froggie_*)\subseteq F_*$ and that $\froggie^+\notin H$.

            If the transition took place according to \cref{rule:cFH}, then $\froggie_*=\froggiecrown_*$ and so clearly $\column(\square)\subseteq F_*$ for every $\square\in I[\froggiecrown_*,\froggie_*]$.
            If the transition took place according to \cref{rule:fFH}, then $\froggiecrown_*=\froggiecrown$, $x=\settled$, thus \cref{cf:two:empty} of $\crown F$ implies $\froggiecrown \notin A$, so $A=\{\froggie\}$.
            Therefore, since $\crown F$ satisfies \cref{cf:one}, we immediately have that $\column(\square)\subseteq F_*$ for every $\square\in I[\froggiecrown_*,\froggie_*]$.

            Now, suppose for the sake of contradiction that $\froggie_*^+\in H_*$ and set $\froggie'=\opp(\froggie_*)$.
            Since $\froggie_*\in H_*$, we know that $\froggie_*^+\neq\froggie'$; phrased differently, $\opp(\froggie')^+\neq\froggie'$.
            This also means that $\opp(\froggie')^+=\froggie_*^+\in H$ since the only difference between $H$ and $H_*$ occurs in $\column(\froggie_*)$.
            However, since $\froggie_*=\froggie^+\notin H$ and $\froggie_*\in F$, this means that $\froggie'\in H$ as well; a contradiction to the fact that $(F,H)$ is a hatted-frog arrangement (\Cref{oppPlus}).

        \item Suppose that $\abs{A_*}=2$.
            Here, \Cref{moveObservations} tells us that $F_*=F$ and $\abs A=2$.
            Furthermore, the transition from $\crown F$ to $\crown F_*$ must have taken place according to \cref{rule:fFH}.
            In particular, $\froggiecrown_*=\froggiecrown$, $A_*=\{\froggiecrown_*,\froggie^+\}$, $x_*=\agitated$ and $\column(\froggie^+)\subseteq F_*$.
            Of course, this means that $F_*\cap \column(\froggiecrown_*)=\varnothing$ and $\column(\square)\subseteq F_*$ for every $\square\in I[\opp(\froggiecrown_*),\froggie^+]-\opp(\froggiecrown_*)$ since $\crown F$ satisfies \cref{cf:two}.

            The argument that $\froggie^{++}\notin H_*$ is identical to the argument used to verify $\froggie^+_* \not\in H_*$ in the preceding item.

        \item Suppose that $x_*=\agitated$ and that $\eb{F_*}(\froggiecrown_*)=\opp(\froggiecrown_*^+)$, yet $F_*\cap\column(\froggiecrown_*^+)\neq\varnothing$.
            Since $x_*=\agitated$, we know by \cref{cf:whereCrown} of $\crown F_*$ (checked above) that $\froggiecrown_*\in A_*$ and so $F_*\neq[2]\times[k]$ (\cref{cf:sum} since $m\leq 2k$); therefore, $\eb{F_*}(\froggiecrown_*)\notin F_*$.
            This implies that $\froggiecrown_*^+\in H_*$.
            We claim also that $\froggiecrown_*\in H_*$.
            Indeed, if $\froggiecrown_*\notin H_*$, then $F_*\cap\column(\froggiecrown_*)=\varnothing$ due to \cref{cf:hatOrEmpty}.
            But then $\eb{F_*}(\froggiecrown_*)=\froggiecrown_*$ and we know that $\froggiecrown_*\neq\opp(\froggiecrown_*^+)$ since $\froggiecrown_*^+\in H_*\subseteq F_*$.

            In particular, this implies that $\froggiecrown_*$ and $\froggiecrown_*^+$ reside in different columns.
            This implies additionally that $\column(\froggiecrown_*)\subseteq F_*$ since we must have $I(\opp(\froggiecrown_*^+),\froggiecrown_*]\subseteq F_*$ (because $\eb{F_*}(\froggiecrown_*)=\opp(\froggiecrown_*^+)$ and $\froggiecrown_*\neq\opp(\froggiecrown_*^+)$)

            We now consider the transition from $\crown F$ to $\crown F_*$.

            If $x=\settled$, then the transition could have only taken place according to \cref{rule:fFC} since $x_*=\agitated$.
            However, in this case, $\opp(\froggiecrown_*),\froggiecrown_*^+\in H$ which contradicts the fact that $(F,H)$ is a hatted-frog arrangement (\Cref{oppPlus}).
            Therefore $x=\agitated$.

            Since $F_*=F+\froggie^+$ and $\column(\froggiecrown_*)\subseteq F_*$, we know that $F\cap\column(\froggiecrown_*)\neq\varnothing$.
            In particular, we cannot have $\abs A=2$ since then $\froggiecrown=\froggiecrown_*$ and so \cref{cf:two:empty} would be violated for $\crown F$.
            Thus, $A=\{\froggie\}$ and $\froggie=\froggiecrown=\froggiecrown_*^-$.
            Since $x_*=\agitated$ as well, we know that $F=F_*$ and so $\eb{F}(\froggiecrown)=\eb{F_*}(\froggiecrown_*)=\opp(\froggiecrown_*^+)$.
            Since $\froggiecrown_*=\froggiecrown^+$, we find that $I(\eb{F}(\froggiecrown),\froggiecrown]$ contains $\opp(\froggiecrown_*)$ but does not contain $\froggiecrown_*$.
            Combining this with the fact that $(F,H)$ aligns with $I(\eb{F}(\froggiecrown),\froggiecrown]$ (\cref{cf:hatOrEmpty}), we must additionally have $\opp(\froggiecrown_*)\in H$.
            But then, both $\opp(\froggiecrown_*)$ and $\froggiecrown_*^+$ belong to $H$; a final contradiction (\Cref{oppPlus}).
    \end{properties}

    \medskip

    Next, we show that $\crown F_*\notin\crownedFrogsStart km$, which will conclude the proof.

    Suppose for the sake of contradiction that $\crown F_*\in\crownedFrogsStart km$.
    Then $A_*\subseteq\column(\froggiecrown_*)$, $x_*=\agitated$ and $F_*$ is disjoint from $\column(\froggiecrown_*)$.
    In particular, $A_*=A-\froggie$ since otherwise we would have $\froggie^+\in A_*\cap F_*$, as per \Cref{moveObservations}.
    Combining these observations, we find that the transition from $\crown F$ to $\crown F_*$ could only have taken place according to \crefOr{rule:fE,rule:fEH}.
    As such, $\froggie^+\in F_*\setminus F$, $\froggiecrown_*=\froggiecrown$ and $x=\agitated$.
    Using rotational symmetry, we may suppose that $\froggiecrown_*=\froggiecrown=(1,c)$.

    Now, according to \cref{cf:two:other}, we know that $\column(\square)\subseteq F$ for every $\square\in I[\opp(\froggiecrown),\froggie]-\opp(\froggiecrown)$.
    Since $\froggie^+\notin F$, this means that $\froggie^+=(2,f)$ for some $f\in[c-1]$.
    Let $t\in[f,c-1]$ be the largest integer for which $(2,f)\in H_*$; note that $t$ exists since $\froggie^+\in H_*$.
    Since $[2]\times[f+1,c-1]\subseteq F_*$, this implies that $(2,t),(1,t+1)\in H_*+\froggiecrown_*$; a contradiction to the fact that $(F_*+\froggiecrown_*,H_*+\froggiecrown_*)$ is a hatted-frog arrangement.
\end{proof}

Now that we know that $\move$ is a function from $\crownedFrogs km\setminus\crownedFrogsEnd km$ to $\crownedFrogs km\setminus\crownedFrogsStart km$, we observe that $\poke{c}$, $\move$ and $\dethrone$ together embody the hatted-frog process.
Indeed, a comparison between \Cref{moveObservations} and the proof of \Cref{tohats} shows that:
\begin{corollary}\label{crown-coupling}
    Fix any $\hat F\in\hattedFrogs km$ and any $c\in[k]$.
    If $\hat F$ has a frog in column $c$, then
    \[
        \hat F c=\bigl(\dethrone\circ\move^h\circ\poke c\bigr)(\hat F),\qquad\text{where }h=\hop(\hat F,c).
    \]
\end{corollary}
Examples of \Cref{crown-coupling} are shown in \Cref{fig:crownedProcess}.
\medskip

At this point, we have ``lifted'' the hatted-frog process to a new ``crowned-frog process'' defined by the map $\move$.
The astute reader will notice that $\move$ is not the only map which would accomplish this feat; in fact, it is far from the simplest such map.
The reason for using this specific transition map is that it has an extra important property:

\begin{theorem}\label{injective}
    For integers $k,m$ with $0\leq m\leq 2k$, $\move\colon\crownedFrogs km\setminus\crownedFrogsEnd km\to\crownedFrogs km\setminus\crownedFrogsStart km$ is an injection.
\end{theorem}

\begin{proof}
    Fix $\crown F=(F,H,A,\froggiecrown,x)\in\crownedFrogs km\setminus\crownedFrogsStart km$ and suppose that $\crown F_1=(F_1,H_1,A_1,\froggiecrown_1,x_1)\in\crownedFrogs km\setminus\crownedFrogsEnd km$ satisfies $\move(\crown F_1)=\crown F$.
    We must show that $\crown F_1$ is unique in this regard.

    Set $\{\froggie_1\}=A_1\setminus A$.

    \textbf{Case 1:} $\abs A=2$.
    Here we have $A=\{\froggiecrown,\froggie\}$ and $\froggie_1=\froggie^-$.
    Furthermore, the transition from $\crown F_1$ to $\crown F$ must have taken place according to \cref{rule:fFH} and so $\crown F_1$ is uniquely determined.
    \medskip

    \textbf{Case 2:} $\abs A=1$ and $x=\settled$.
    In this case, we must have $A_1=\{\froggie_1\}$, $A=\{\froggie_1^+\}$ and $\froggie_1^+\in F$.
    In particular, $\froggie_1$ and $F$ are uniquely determined.
    Recall that \Cref{moveObservations} implies that $\froggiecrown\in\{\froggiecrown_1,\froggie_1^+\}$.

    If $\froggiecrown=\froggie_1^+$, then the transition from $\crown F_1$ to $\crown F$ could only have taken place according to \cref{rule:cFH}.
    If $\froggiecrown\neq\froggie_1^+$, then $\froggiecrown_1=\froggiecrown$ and transition from $\crown F_1$ to $\crown F$ could only have taken place according to \cref{rule:fFH}.
    Therefore, since $\froggie_1$ is uniquely determined, $\crown F_1$ is uniquely determined based on whether or not $\froggiecrown=\froggie_1^+$ (i.e.\ whether or not $\froggiecrown\in A$).
    \medskip

    \textbf{Case 3:} $\abs A=1$ and $F\cap\column(\froggiecrown)=\varnothing$.
    In this case, the transition from $\crown F_1$ to $\crown F$ must have taken place according to \crefOr{rule:fE,rule:fEH} and $\froggie_1^+\neq\froggiecrown$.
    In particular, we know that $x_1=x=\agitated$, $\froggiecrown_1=\froggiecrown$ and $\{\froggie_1^+\}=F\setminus F_1$.
    Therefore, if we can show that $\froggie_1$ is uniquely determined, then we have shown that $\crown F_1$ is uniquely determined.

    Using rotational symmetry, we may suppose that $\froggiecrown=(1,c)$ and so $\froggie_1^+=(r,f)$ for some $f\in[c-1]$ due to \cref{cf:two}.
    If $\froggie_1=(r',f')$, then \cref{cf:two} implies additionally that $[2]\times[f',c-1]\subseteq F_1$.
    Since $(r,f)\notin F_1$, this means that $f'>f$ and so $r=2$ and $f'=f+1$.

    Now, let $t\in[c]$ be the largest integer for which either $(1,t)\notin F+(1,c)$ or $(2,t-1)\notin H$.
    Note that $t$ exists since $(2,0)\notin H$ vacuously.
    Additionally, note that we have $t=c$ if and only if $(2,c-1)\notin H$.
    If we can show that $f=t$, then we will have shown that $\froggie_1$ is uniquely determined as desired.

    As discussed above, we know that $[2]\times[f+1,c-1]\subseteq F$ and know also that $(F,H)$ aligns with $I((2,c),(2,f)]$, implying that $\{2\}\times[f,c-1]\subseteq H$.
    Therefore, we must have $f\geq t$.
    If $f>t$, then it must be the case that $(1,f)\in F$ and $(2,f-1)\in H$.
    However, since $(2,f)\notin F_1$ this would mean that $(2,f-1),(1,f)\in H_1$, contradicting the fact that $(F_1,H_1)$ is a hatted-frog arrangement.
    \medskip

    \textbf{Case 4:} $\abs A=1$, $x=\agitated$ and $F\cap\column(\froggiecrown)\neq\varnothing$.
    Therefore, we must have $A=\{\froggiecrown\}$ and $\froggiecrown\in H$.
    We claim furthermore, that $\froggiecrown=\froggie_1^+$.
    Indeed, if this is not the case, then $\froggiecrown_1=\froggiecrown$ and $A_1=\{\froggiecrown,\froggie\}$.
    This, however, would contradict \cref{cf:two} since we would have $\froggiecrown_1\in F_1\cap\column(\froggiecrown_1)$

    Therefore, either $A_1=\{\froggie_1\}$ or $A_1=\{\froggiecrown,\froggie_1\}$ where $\froggiecrown=\froggie_1^+$.
    In particular, $\froggie_1=\froggiecrown^-$ is uniquely determined.
    Note that if $\abs{A_1}=2$, then we must have $\froggiecrown_1=\froggiecrown$ and $F_1\cap\column(\froggiecrown)=\varnothing$.

    Consider $\eb{F}(\froggiecrown)$.
    Since $F_1\subseteq F\neq[2]\times[k]$ (\cref{cf:sum}), we know that $\eb{F}(\froggie_1)=\eb{F}(\froggiecrown)$ due to \cref{cf:one}.
    \medskip

    \textbf{Case 4a:} $\eb{F}(\froggiecrown)=\opp(\froggiecrown)$.
    In this case, we claim that $\abs {A_1}=2$.
    Indeed, suppose that $A_1=\{\froggie_1\}$.
    If $x_1=\agitated$, then $\froggie_1=\froggiecrown_1$ and we must have $\froggiecrown=\froggiecrown_1^+\in H_1$ or else $x=\settled$.
    However, since $\eb{F_1}(\froggiecrown_1)=\opp(\froggiecrown_1^+)$ here, this contradicts \cref{cf:agitated}.
    On the other hand, if $x_1=\settled$, then the transition from $\crown F_1$ to $\crown F$ must have taken place according to \cref{rule:fFC}; however this contradicts the fact that $\opp(\froggiecrown)=\eb{F}(\froggiecrown)\notin F_1$.

    Therefore, $\abs {A_1} =2$ and so also $F_1\cap\column(\froggiecrown)=\varnothing$.
    Because also $\froggiecrown=\froggie_1^+$, the transition from $\crown F_1$ to $\crown F$ must have taken place according to \cref{rule:fE} and so $\crown F_1$ is uniquely determined.
    \medskip

    \textbf{Case 4b:} $\eb{F}(\froggiecrown)\neq\opp(\froggiecrown)$.
    In this case, we must have $\abs {A_1}=1$.
    Indeed, if $\abs {A_1}=2$, then \cref{cf:hatOrEmpty,cf:two} would imply that $\eb{F}(\froggiecrown)=\opp(\froggiecrown)$.

    Since $\abs{A_1}=1$, the transition from $\crown F_1$ to $\crown F$ could only have taken place according to either \cref{rule:cH,rule:fFC}.
    Note that if $x_1=\agitated$, then $\froggiecrown_1=\froggie_1$ and that if $x_1=\settled$, then $\froggiecrown_1=\opp(\froggiecrown)$.
    Thus, if we can show that $x_1$ is uniquely determined, then we have shown that $\crown F_1$ is uniquely determined as well.

    Recall that $\eb{F}(\froggie_1)=\eb{F}(\froggiecrown)$.
    We next claim that also $\eb{F}(\froggiecrown_1)=\eb{F}(\froggiecrown)$.
    Indeed, if $x_1=\agitated$, then $\froggiecrown_1=\froggie_1$ and so $\eb{F}(\froggiecrown_1)=\eb{F}(\froggiecrown)$
    Otherwise, if $x_1=\settled$, then \cref{cf:one} tells us that $I[\froggiecrown_1,\froggie_1]\subseteq F$ and so we still have that $\eb{F}(\froggiecrown_1)=\eb{F}(\froggiecrown)$.

    Using rotational symmetry, we may suppose that $\froggiecrown=(1,c)$.
    Set $\eb{F}(\froggiecrown)=(r,f)$; since $F=F_1\neq[2]\times[k]$, we know that $(r,f)\notin F$.
    In particular, $f\neq c$.
    We break into two cases:
    \medskip

    \textbf{Case 4bi:} $f<c$.
    We claim that $x_1=\agitated$.
    Indeed, if $x_1=\settled$, then $\froggiecrown_1=(2,c)$ and so $(r,f)\in I[\froggiecrown_1,\froggie_1]\subseteq F_1=F$; a contradiction.
    \medskip

    \textbf{Case 4bii:} $f>c$.
    In this case, we know that $(2,c)\in I((r,f),\froggie_1]$ whereas $(1,c)\notin I((r,f),\froggie_1]$.
    Since $(F_1,H_1)$ aligns with $I((r,f),\froggie_1]$, this implies that $(2,c)\in H$ and so $x_1=\settled$ since otherwise $(1,c)\in H$.
    \medskip

    \textbf{Case 5:} $\abs A=0$.
    In this case, we must have $A_1=\{\froggie_1\}$ and $\{\froggie_1^+\}=F\setminus F_1$.
    In particular, the transition from $\crown F_1$ to $\crown F$ took place according to one of \crefOr{rule:cE,rule:cEH,rule:fE,rule:fEH}.
    As such, if we can show that $\froggie_1$ is uniquely determined, then we will have shown that $\crown F_1$ is uniquely determined.
    \medskip

    \textbf{Case 5a:} $F=[2]\times[k]$.
    Using rotational symmetry, we may suppose that $\froggiecrown=(1,c)$ for some $c\in[k]$.
    Since $\froggiecrown\in H$ and $F$ contains every square of the grid, we know that there is some $\ell\in[k]$ for which $(1,f)\in H$ whenever $f\leq\ell$ and $(2,f)\in H$ whenever $f>\ell$.
    Of course, $c\leq\ell$.
    We claim that $\froggie_1^+=(1,\ell)$ which will conclude this case.

    To see why, because $\froggie_1^+\in H$, we have either $\froggie_1^+=(1,f)$ for some $f\leq\ell$ or $\froggie_1^+=(2,f)$ for some $f>\ell$.
    In the former case, we would have $(2,f)\in H_1$ and so we must also have $(2,\ell)\in H_1$ implying that $f=\ell$.
    In the latter case, we would have $\froggie_1\neq\froggiecrown_1=(1,c)$ yet $(2,f)\notin F_1$ which contradicts \cref{cf:one}.
    \medskip

    \textbf{Case 5b:} $F\neq[2]\times[k]$.
    Using rotational symmetry, we may suppose that $\froggiecrown=(1,c)$.
    Let $f\in[c,k]$ be the smallest integer such that either $(2,f)\notin F$ or $(1,f)^+\notin H$.
    By construction, $\{1\}\times[c,f]\subseteq H$ and $[2]\times[c,f-1]\subseteq F$.
    We claim that $\froggie_1^+=(1,f)$, which will conclude this case.

    In this situation we know that $\eb{F}(\froggie_1^+)=\eb{F}(\froggiecrown)$ and that $(F,H)$ aligns with $I(\eb{F}(\froggiecrown),\froggie_1^+]$.
    Since $\froggie_1^+\notin F_1$ we must have $\froggie_1^+=(1,\ell)$ for some $\ell\in[c,k]$ or else \cref{cf:one} is violated for $\crown F_1$.
    Further using \cref{cf:one}, $\{1\}\times[c,\ell]\subseteq H$ and that $[2]\times[c,\ell-1]\subseteq F$.
    In particular, $\ell\geq f$.
    Certainly if $(2,\ell)\notin F$, then $\ell=f$.
    Otherwise, $(2,\ell)\in F$ and hence $(2,\ell)\in H_1$.
    Since $H_1$ and $H$ differ only in column $\ell$, this would mean that $(1,\ell+1)\in H_1$ in the case that $\ell<f$, yielding a contradiction.
\end{proof}

The injectivity of $\move$ can now be easily bootstrapped up to full bijectivity.

\begin{corollary}\label{moveBijection}
    Fix integers $k,m$ with $0\leq m\leq 2k$.
    The map $\move\colon\crownedFrogs km\setminus\crownedFrogsEnd km\to\crownedFrogs km\setminus\crownedFrogsStart km$ is a bijection.
    Moreover, for any $\crown F\in\crownedFrogs km$, there is a \emph{unique} pair $(\crown S,h)\in\crownedFrogsStart km\times\Z_{\geq 0}$ for which $\crown F=\move^h(\crown S)$.
\end{corollary}
\begin{proof}
    We begin by observing that for any $\crown F\in\crownedFrogs km$, there is a unique $t\in\Z_{\geq 0}$ for which $\move^t(\crown F)\in\crownedFrogsEnd km$.
    Indeed, the set of agitated frogs is weakly-decreasing in size and strictly decreases when a frog hops to an empty square.
    Since $m\leq 2k$, there is always an empty square if there is an agitated frog and so the observation follows.

    Now, observe that for $(F,H)\in\hattedFrogs km$ and $c\in[k]$, $F$ contains an element of column $c$ if and only if $H$ does as well.
    As such, \Cref{injective-poke,ending-select} together imply that
    \[
        \abs{\crownedFrogsStart km}=\sum_{(F,H)\,\in\,\scalebox{0.65}{$\hattedFrogs km$}}\abs{H}=\abs{\crownedFrogsEnd km}.
    \]
    Together with the observation in the previous paragraph and \Cref{injective}, this establishes the claim.
\end{proof}

The proof of \Cref{regular} is now almost immediate.

\begin{proof}[Proof of \Cref{regular}]
    Fix any $\hat F=(F,H)\in\hattedFrogs km$; we must show that there are exactly $\abs\bet$ many pairs $(\hat S,a)\in\hattedFrogs km\times\bet$ for which $\hat Sa=\hat F$.

    Note that $k-\abs H$ is precisely the number of empty columns in $\hat F$.
    Of course, if $c\in\bet$ is either the index of an empty column in $\hat F$ or not the index of a column (i.e.\ $c\notin[k]$), then $\hat Fc=\hat F$.
    This accounts for exactly $\abs\bet-\abs H$ many in-edges to $\hat F$.

    Therefore, we must show that there are exactly $\abs H$ many pairs $(\hat S,c)\in\hattedFrogs km\times[k]$ for which both $\hat Sc=\hat F$ and $\hat S$ has a frog in column $c$.
    Since $\abs{\dethrone^{-1}(\hat F)}=\abs H$ (\Cref{ending-select}), this is an immediate consequence of \Cref{crown-coupling,moveBijection}.
\end{proof}

\section{Computation of the speeds}\label{sec:speeds}

Now that we know that the $m$-hatted-frog dynamics admits a uniform stationary distribution, we can exploit this fact to compute the speeds of the frogs in the original frog dynamics.

\begin{theorem}\label{cumulative}
    Set $W=12\cdots kk\cdots 21$ and let $\bet\supseteq[k]$ be any alphabet.
    For each $m\in[2k]$, the cumulative speed of the nastiest $m$ frogs in the $(W,\bet)$ frog dynamics is
    \[
        \sum_{i=1}^m s_i={2k\cdot\sum_{j\geq 0}{2k-2j\choose m-1-2j}\over\abs\bet\cdot\sum_{j\geq 0}{2k-2j\choose m-2j}}.
    \]
\end{theorem}

Together with \Cref{frogTheorem}, this proves \Cref{lcsResult}.
The proof of \Cref{cumulative} encapsulates the whole of this section.
\medskip

We begin by detailing an important bijection, which we derive by applying our knowledge of the crowned frogs from the previous section.
\begin{lemma}\label{halfway}
    Fix integers $k,m$ with $0\leq m\leq 2k$.
    For every $\square\in[2]\times[k]$, there is a bijection
    \[
        \phi_\square\colon
        \bigl\{(F,H)\in\hattedFrogs km: \square\in H,\ I[\opp(\square),\square]\not\subseteq F\bigr\}\to
        \bigl\{(F,H)\in\hattedFrogs k{m-1}: F\cap\column(\square)=\varnothing\bigr\}.
    \]
\end{lemma}
Informally, the bijection is accomplished by poking the column containing $\square$ and following the hatted-frog process until exactly one frog is agitated.
\begin{proof}
    Fix $\hat F=(F,H)\in\hattedFrogs km$ with $\square\in H$ and $I[\opp(\square),\square]\not\subseteq F$.
    Suppose that $\square=(r,c)$.
    For $t\in[0,\hop(\hat F,c)]$, set $(F_t,H_t,A_t,\froggiecrown_t,x_t)=\move^t(\poke{c}(\hat F))$, which is reasonable due to \Cref{crown-coupling}.
    Let $T\in[0,\hop(\hat F,c)]$ be the smallest integer for which $\abs{A_T}=1$ and define $\phi_\square(\hat F)\eqdef(F_T,H_T)$.
    The fact that $(F_T,H_T)\in\hattedFrogs k{m-1}$ and $F_T\cap\column(\square)=\varnothing$ follow from \Cref{well-defined} and the fact that $\abs{F_T}+\abs{A_T}=m$.

    To see why $\phi_\square$ is bijective, fix any $(F',H')\in\hattedFrogs k{m-1}$ with $F'\cap\column(\square)=\varnothing$.
    Then $\crown F=(F',H',\{\square\},\square,\agitated)$ is a member of $\crownedFrogs km$.
    By \Cref{injective-poke} and \Cref{moveBijection}, there is a unique $(\hat F,c,t)\in\hattedFrogs km\times[k]\times\Z_{\geq 0}$ for which $\crown F=\move^t(\poke{c}(\hat F))$.
    If $\hat F=(F,H)$, then it is routine to check that $\square\in H$ and $I[\opp(\square),\square]\not\subseteq F$.
\end{proof}

The above $\phi_\square$ is actually just a piece of a more far-reaching bijection which we now detail.

Fix a hatted-frog arrangement $\hat F\in\hattedFrogs km$ and fix an index $c\in[k]$.
Define $\Hop(\hat F,c)\subseteq[2]\times[k]$ to be the set of frogs that hopped when poking column $c$ in the arrangement $\hat F$ in the hatted-frog process.
Formally, if $\hat F=(F,H)$, then
\[
    \Hop(\hat F,c)\eqdef\bigl\{\froggie\in[2]\times[k]:I[(r,c),\froggie]\subseteq F\text{ for some }r\in[2]\bigr\}.
\]
By definition, $\hop(\hat F,c)=\abs{\Hop(\hat F,c)}$.

Expanding upon this, for non-negative integers $k,m$, define
\[
    \Omega_{k,m}\eqdef\bigl\{(\hat F,c,\froggie):\hat F\in\hattedFrogs km,\ c\in[k],\ \froggie\in\Hop(\hat F,c)\bigr\},
\]
and observe that
\begin{equation}\label{eqn:speedSum}
    \abs{\Omega_{k,m}}=\sum_{c\in[k]}\ \sum_{\scalebox{0.65}{$\hat F\!\in\!\hattedFrogs km$}}\hop(\hat F,c).
\end{equation}
Due to \Cref{uniform}, computing $\abs{\Omega_{k,m}}$ is essentially equivalent to proving \Cref{cumulative}.

We next define two maps
\begin{align*}
    \Phi_{k,m} &\colon \Omega_{k,m}\to\hattedFrogs k{m-1}\times([2]\times[k]),\\
    \Psi_{k,m} &\colon \hattedFrogs k{m-1}\times([2]\times[k])\to\Omega_{k,m},
\end{align*}
which will be used in our remaining arguments.
\medskip

Fix $(\hat F,c,\froggie)\in\Omega_{k,m}$ and suppose that $\hat F=(F,H)$.
It must be the case that $F$ contains some element of column $c$, otherwise $\froggie$ would not exist.
Thus, set $\{\froggiehat\}=H\cap([2]\times\{c\})$.
\begin{enumerate}
    \item If $I[\opp(\froggiehat),\froggie]\subseteq F$, then
        \[
            \Phi_{k,m}(\hat F,c,\froggie)\eqdef\bigl((F-\opp(\froggiehat),H),\froggie\bigr).
        \]
    \item Otherwise, $I[\opp(\froggiehat),\froggie]\not\subseteq F$, so define
        \[
            \Phi_{k,m}(\hat F,c,\froggie)\eqdef\bigl(\phi_{\scalebox{0.75}{$\froggiehat$}}(\hat F),\froggie\bigr),
        \]
        where $\phi_{\scalebox{0.75}{$\froggiehat$}}$ is the bijection from \Cref{halfway}.
\end{enumerate}
Due to \Cref{halfway}, it is easy to see that $\Phi_{k,m}\colon\Omega_{k,m}\to\hattedFrogs k{m-1}\times([2]\times[k])$ is well-defined.
\medskip

Next, fix $(\hat F,\square)\in\hattedFrogs k{m-1}\times([2]\times[k])$ and suppose that $\hat F=(F,H)$.
Note that $\eb{F}(\square)\notin F$ since $m-1<2k$.
Suppose that $\eb{F}(\square)=(r,c)\in[2]\times[k]$.
\begin{enumerate}
    \item If $\opp(\eb{F}(\square))\in F$ (and so $\opp(\eb{F}(\square))\in H$), then
        \[
            \Psi_{k,m}(\hat F,\square)\eqdef\bigl((F+\eb{F}(\square),H),c,\square\bigr).
        \]
    \item Otherwise, $F\cap\column(\eb{F}(\square))=\varnothing$, so define
        \[
            \Psi_{k,m}(\hat F,\square)\eqdef\bigl(\phi_{\eb{F}(\square)}^{-1}(\hat F),c,\square\bigr),
        \]
        where $\phi_{\eb{F}(\square)}$ is the bijection from \Cref{halfway}.
\end{enumerate}
Due to \Cref{halfway}, it is again easy to see that $\Psi_{k,m}\colon\hattedFrogs k{m-1}\times([2]\times[k])\to\Omega_{k,m}$ is well-defined.
\medskip

Observe that for any $F\subseteq[2]\times[k]$ and any $\square_1,\square_2\in[2]\times[k]$, we have $I[\square_1,\square_2]\subseteq F$ if and only if $\eb{F-\square_1}(\square_2)=\square_1$.
This observation along with \Cref{halfway} implies:
\begin{theorem}\label{speedBijection}
    For any integers $k,m$ with $0\leq m\leq 2k$, the maps $\Phi_{k,m}$ and $\Psi_{k,m}$ are inverses.
\end{theorem}

Immediately, \Cref{speedBijection} implies that $\Phi_{k,m}\colon\Omega_{k,m}\to\hattedFrogs k{m-1}\times([2]\times[k])$ is a bijection whenever $0\leq m\leq 2k$.
Combining this fact with \cref{eqn:speedSum}, we immediately have:
\begin{theorem}\label{speeds}
    For any integers $k,m$ with $0\leq m\leq 2k$,
    \[
        \sum_{c\in[k]}\ \sum_{\scalebox{0.65}{$\hat F\!\in\!\hattedFrogs km$}}\hop(\hat F,c)=2k\cdot\abs{\hattedFrogs k{m-1}}.
    \]
\end{theorem}

By taking appropriate ``slices'' of $\Phi_{k,m}$ and $\Psi_{k,m}$, \Cref{speedBijection} has another straightforward implication.
\begin{lemma}\label{corner}
    For any integers $k,m$ with $0\leq m\leq 2k-1$,
    \[
        \bigl\lvert\bigl\{(F,H)\in\hattedFrogs km:(1,k)\in F\bigr\}\bigr\rvert=\bigl\lvert\bigl\{(F,H)\in\hattedFrogs k{m-1}:(1,k)\notin F\bigr\}\bigr\rvert
    \]
\end{lemma}
\begin{proof}
    Define
    \begin{align*}
        \mathcal A &\eqdef\bigl\{\bigl((F,H),k,(1,k)\bigr):(F,H)\in\hattedFrogs km,\ (1,k)\in F\bigr\},\quad\text{and}\\
        \mathcal B &\eqdef\bigl\{\bigl((F,H),(1,k)\bigr):(F,H)\in\hattedFrogs k{m-1},\ (1,k)\notin F\bigr\}.
    \end{align*}
    To prove the lemma, it suffices to show that $\abs{\mathcal A}=\abs{\mathcal B}$.
    Observe that $\mathcal A\subseteq\Omega_{k,m}$ and that $\mathcal B\subseteq\hattedFrogs k{m-1}\times([2]\times[k])$.

    Now, fix any $(F,H)\in\hattedFrogs km$ with $(1,k)\in F$.
    Since $m<2k$, we know that $I[(2,k),(1,k)]\not\subseteq F$ and therefore $\Phi_{k,m}\bigl((F,H),k,(1,k)\bigr)\in\mathcal B$.

    On the other hand, fix any $(F,H)\in\hattedFrogs k{m-1}$ with $(1,k)\notin F$.
    Then $\eb{F}(1,k)=(1,k)$ and so $\Psi_{k,m}\bigl((F,H),(1,k)\bigr)\in\mathcal A$.

    Therefore, $\Phi_{k,m}\colon\mathcal A\to\mathcal B$ and $\Psi_{k,m}\colon\mathcal B\to\mathcal A$, and so the claim follows from \Cref{speedBijection}.
\end{proof}

We can now compute the total number of hatted-frog arrangements.

\begin{theorem}\label{size}
    For any non-negative integers $k,m$,
    \[
        \abs{\hattedFrogs km}=\sum_{i\geq 0}{2k-2i\choose m-2i}.
    \]
\end{theorem}
Before proceeding with the proof, we note that these same numbers have appeared in various other contexts; see \oeis{A035317}, \oeis{A059259} and \oeis{A059260}.
\begin{proof}
    For non-negative integers $n,m$, define
    \[
        f(n,m)\eqdef\begin{cases}
            \abs{\hattedFrogs{k}{m}}, &\text{if }n=2k,\\
            \bigl\lvert\bigl\{(F,H)\in\hattedFrogs {k+1}m:(1,k+1)\notin F\bigr\}\bigr\rvert, &\text{if }n=2k+1.
        \end{cases}
    \]
    We will prove that
    \begin{equation}\label{eqn:identity}
        f(n,m)=\sum_{i\geq 0}{n-2i\choose m-2i},
    \end{equation}
    which will establish the claim since $f(2k,m)=\abs{\hattedFrogs km}$ for any non-negative integer $k$.

    It is not difficult to see that $f(n,0)=1$ and $f(n,n)=\lfloor n/2\rfloor+1$ for any non-negative $n$.
    We will show that $f$ additionally satisfies the recurrence
    \begin{equation}\label{eqn:recurrence}
        f(n,m)=f(n-1,m)+f(n-1,m-1),\quad\text{whenever }0<m<n.
    \end{equation}
    Since the same boundary conditions and recurrence is true of the numbers $\sum_{i\geq 0}{n-2i\choose m-2i}$, this will establish \cref{eqn:identity} and conclude the proof.
    \medskip

    If $n$ is odd, say $n=2k+1$, then clearly
    \[
        \bigl\{(F,H)\in\hattedFrogs{k+1}m:(1,k+1)\notin F\bigr\} = \hattedFrogs km\sqcup\bigl\{\bigl(F+(2,k+1),H+(2,k+1)\bigr):(F,H)\in\hattedFrogs k{m-1}\bigr\},
    \]
    which establishes \cref{eqn:recurrence} in this case.

    Next, if $n$ is even, say $n=2k$, then
    \[
        \hattedFrogs km = \bigl\{(F,H)\in\hattedFrogs km:(1,k)\notin F\bigr\}\sqcup\bigl\{(F,H)\in\hattedFrogs km:(1,k)\in F\bigr\}.
    \]
    The former set has size $f(n-1,m)$ by definition, and the latter set has size $f(n-1,m-1)$ due to \Cref{corner} and so we have established \cref{eqn:recurrence}.
\end{proof}

We finally end our journey by combining all of the prior results in this paper to deduce the speeds of the frogs in the frog dynamics associated with $W=12\cdots kk\cdots 21$:

\begin{proof}[Proof of \Cref{cumulative}]
    Recall that $W=12\cdots kk\cdots 21$, that $\bet\supseteq[k]$ is some alphabet and that $m\in[2k]$.
    Let $\pi$ denote the stationary distribution of the $(W,\bet)$ $m$-blind-frog dynamics; \Cref{speeds_v_hops} tells us that
    \[
        \sum_{i=1}^m s_i=\E_{a\sim\bet}\ \E_{F\sim\pi} \hop(F,a).
    \]
    Now, due to the coupling in \Cref{coupling} and the fact that the $(k,\bet)$ $m$-hatted-frog dynamics admits a uniform stationary (\Cref{uniform}), we know that
    \begin{align*}
        \E_{a\sim\bet}\ \E_{F\sim\pi_m} \hop(F,a) &= \E_{a\sim\bet}\ \E_{\scalebox{0.65}{$\hat F\!\sim\!\hattedFrogs km$}}\hop(\hat F,a) = {1\over\abs\bet\cdot\abs{\hattedFrogs km}}\sum_{a\in\bet}\ \sum_{\scalebox{0.65}{$\hat F\!\in\!\hattedFrogs km$}}\hop(\hat F,a)\\
                                                  &= {1\over\abs\bet\cdot\abs{\hattedFrogs km}}\sum_{a\in[k]}\ \sum_{\scalebox{0.65}{$\hat F\!\in\!\hattedFrogs km$}}\hop(\hat F,a)= {2k\cdot\sum_{j\geq 0}{2k-2j\choose m-1-2j}\over\abs\bet\cdot\sum_{j\geq 0}{2k-2j\choose m-2j}},
    \end{align*}
    where the final equality follows from \Cref{speeds,size}.
\end{proof}

\paragraph{Acknowledgments.}
The frog symbol is from \emph{Froggy Font} (\url{https://www.dafont.com/froggy.font}) by Vladimir Nikolic.
The lily pad symbol is based on a drawing by FrauBieneMaja (\url{https://pixabay.com/vectors/water-lily-lake-water-pond-blossom-4177686/}).

\paragraph{Competing interests.}
The authors declare none.

\bibliographystyle{abbrv}
\bibliography{references}
\end{document}